\documentclass[12pt]{article}
\usepackage{amsmath,enumerate,amsfonts,amssymb}
\usepackage[usenames,dvipsnames]{color}

\usepackage{ntheorem}
\usepackage{graphicx}
\usepackage{graphics, float}

\def\ZZ{{\mathbb Z}}
\def\RR{{\mathbb R}}

\def\mcF{{\mycal F}}

\def\mcP{{\mycal P}}

\newcommand{\mcI}{\mathcal{I}}

\newcommand{\Qnr}{Q^t_{NR}}
\newcommand{\Qnrz}{Q^v_{NR}}

\def\f{{\varphi}}

\newtheorem{theorem} {\sc  Theorem\rm} [section]

\newtheorem{lemma} [theorem] {\sc  Lemma\rm}

\newtheorem{definition}[theorem]{\sc  Definition\rm}

\theorembodyfont{\upshape}
\newtheorem{remark}[theorem]{\sc  Remark\rm}

\def\bproof{\noindent{\bf Proof.\;}}
\def\eproof{\hfill$\square$\medskip}

\newcounter{marnote}

\DeclareFontFamily{OT1}{rsfs}{}
\DeclareFontShape{OT1}{rsfs}{m}{n}{ <-7> rsfs5 <7-10> rsfs7 <10-> rsfs10}{}
\DeclareMathAlphabet{\mycal}{OT1}{rsfs}{m}{n}

\newcommand{\defeq}{\stackrel{\scriptscriptstyle \text{def}}{=}}

\def\tr{{\rm tr}}

\def\mcS{{\mycal{S}}}

\def\mcE{{\mycal E}}

\def\be{\begin{equation}}
\def\ee{\end{equation}}
\def\bea#1\eea{\begin{align}#1\end{align}}
\def\non{\nonumber}
\newcommand{\R}{\mathbb{R}}

\definecolor{orange}{RGB}{255,127,0}

\numberwithin{equation}{section}
\allowdisplaybreaks

\begin{document}

\title{Liquid crystal defects in the Landau-de Gennes theory in two dimensions -- beyond the one-constant approximation}

\author{ Georgy Kitavtsev\thanks{School of Mathematics, University of Bristol, University Walk,
Bristol, BS8 1TW, UK} \hspace{-1ex}, Jonathan M Robbins\footnotemark[1] \hspace{-0.7ex},\\ Valeriy Slastikov\footnotemark[1] \hspace{-0.7ex},    Arghir Zarnescu\thanks{ IKERBASQUE, Basque Foundation for Science, Maria Diaz de Haro 3,
48013, Bilbao, Bizkaia, Spain}\,\,\thanks{BCAM, Basque Center for Applied Mathematics, Mazarredo 14, E48009 Bilbao, Basque Country, Spain}\,\,\thanks{ ``Simion Stoilow" Institute of the Romanian Academy, 21 Calea Grivi\c{t}ei,  010702 Bucharest, Romania}}

\maketitle
\begin{abstract}
We consider the two-dimensional Landau-de Gennes energy with several elastic constants, subject to general $k$-radially symmetric boundary conditions. We show that for generic elastic constants the critical points consistent with the symmetry of the boundary conditions exist only in the case $k=2$. In this case we identify three types of radial profiles: with two, three of full five components and numerically investigate their minimality and stability depending on suitable parameters. 

We also numerically study the stability properties of the critical points of the Landau-de Gennes energy and capture the intricate dependence of various qualitative features of these solutions on the elastic constants and the physical regimes of the liquid crystal system.
\end{abstract}
\section{Introduction }

Nematic liquid crystals are amongst the simplest complex materials. They are liquid solutions  whose microstructure is determined by the  local orientational order of their  rod-like solute particles. This microstructure is responsible for the anisotropic optical properties on which LCD  technologies are based.  
Defects in nematic liquid crystals are localised regions 
where the microstructure varies rapidly, creating striking  patterns in transmitted polarised light.  
An account of defects is a basic challenge for any theoretical description of nematic liquid crystals \cite{chandra, dg}.

In the continuum setting, the material is described by functions taking values in an {\it order-parameter space}, which describes the microstructure. The choice of order-parameter space varies between models.  
In the simplest model, the Oseen-Frank theory, the order-parameter  space is  the unit sphere $\mathbb{S}^2$ with the unit-length vector, $n(x)$, representing the mean local orientation of the  molecules at any point $x$. Equilibrium configurations correspond to minimisers of an energy described  by a local energy density $\sigma(n,\nabla n)$.  Defects correspond to discontinuities in $n$, whose existence may be enforced, for topological reasons, by  boundary conditions.  Point defects have finite energy in three dimensions, but point defects in two dimensions (or line defects in three dimensions)   have infinite energy; this is one of the limitations of the Oseen-Frank theory.  

A more refined (and thus more complicated) model is provided by the Landau-de Gennes theory.  Here, the order-parameter space is the set of traceless symmetric $3\times 3$ matrices, referred to in this context as  $Q$-tensors, which describe the second moments of the orientational probability distribution for the solute particles 
 (see, eg,  \cite{mottram2014introduction}).  Equilibrium configurations $Q(x)$ are minimisers of an energy described by a local  energy density $f(Q,\nabla Q)$.  The Oseen-Frank theory   emerges  in the physically relevant regime where $Q$ has (nearly) everywhere a (nearly) doubly degenerate eigenvalue of fixed magnitude; $Q$ is then (nearly) determined by the orientation of its unique nondegenerate eigenvector, which corresponds to the Oseen-Frank director $n$.  Defects correspond to regions where $Q$ is far from having a doubly degenerate eigenvalue.  The energy of these regions is finite, and  the Landau-de Gennes theory provides a resolution of the discontinuities in the Oseen-Frank model into smooth localised spatial profiles.  

The analytic description of defects in the Landau-de Gennes theory is a formidable mathematical challenge.  Nevertheless, there has been substantial progress in recent years  in analysing the existence, structure,  stability and instability of certain defects in three dimensions \cite{Ball-Maj, Ball-Zar,INSZ_CRAS, INSZ3, lamy2015uniaxial, Lin-Liu, Ma-Za} and two dimensions \cite{baumanphilips, canevari, DRSZ, onsager-spatial_variation, Fat-Slas, GolovatyMontero, disclinationsnumerics, INSZ2d1, INSZ2d2}.
These studies have focussed almost exclusively on a simplified version of the Landau-de Gennes theory, in which the elastic (ie, $\nabla Q$-dependent) contribution to the energy density consists solely of the Dirichlet energy  -- this is the so-called one-constant approximation.  
It is well known that symmetry considerations allow for additional elastic energy terms.  

Here we undertake an exploration of defects in the Landau-de Gennes theory  in two-dimensional domains when the one-constant restriction is relaxed.  Mathematically, this introduces new challenges: 
the Euler-Lagrange equations become a genuine  system of PDEs, and many basic techniques valid for scalar equations fail (though some indeed survive in the two-dimensional setting, see for instance \cite{baumanphilips, bauman2015regularity, cubicinstability}).

\subsection{Isotropy in the Oseen-Frank and Landau-de Gennes theories}

\label{sect:several}

In the absence of a microscopic derivation, the energy  in a continuum theory is generally chosen to be i) compatible with underlying symmetries, and ii) as simple as possible while providing a good description of observed phenomena.    For nematic liquid crystals, the energy  is generally expressed in terms of  a local  density which, in the absence of external fields, depends only on the order parameter and its first spatial derivatives (and is thereby compatible with translational symmetry in physical space), and contains terms up to quadratic order in the first derivatives.  (For problems with free boundary conditions, a surface energy  may also be introduced.)

Symmetry under rotations and reflections, which generate the full orthogonal group $O(3)$, imposes additional restrictions on the energy density.    The orthogonal group acts  on physical space as $x\mapsto Px$,  where $P \in O(3)$, 
and separately on order-parameter space according to the theory.  
For the Oseen-Frank theory, the $O(3)$-action is given by $n \mapsto Pn$,  and for the  Landau-de-Gennes theory, by $Q \mapsto P Q P^t$.

If the energy density is required to be \emph {separately} invariant under the $O(3)$-actions on physical space and order-parameter  space, we obtain its simplest and most restricted  form.
For the Oseen-Frank theory,  this requirement implies that the  energy density $\sigma(n,\nabla n)$ is   proportional to the Dirichlet energy 
$n_{i,j} n_{i,j}$.\footnote{Here and in what follows, we use the summation convention for repeated indices; furthermore we denote $n_{i,j} := \frac{\partial n_i}{\partial x_j}$, $\partial_k Q_{ij}=Q_{ij,k}:=\frac{\partial Q_{ij}} {\partial x_k}$.}
For the Landau-de Gennes theory,  
we first decompose the energy density $\psi(Q,\nabla Q)$ into a sum of a bulk term, which depends only on $Q$, and an elastic term, which depends on $Q$ and $\nabla Q$, as follows 
\cite{BallLect}: 
$$\psi(Q,\nabla Q)=\underbrace{\psi(Q,0)}_{:=f_B(Q)\textrm{ bulk part}}+ \ \  \ \underbrace{\psi(Q,\nabla Q)-\psi(Q,0)}_{:=f_E(Q,\nabla Q)\textrm{ elastic part }}.$$ 
Invariance under the separate $O(3)$-actions  implies that the elastic energy density $f_E$ is proportional to the Dirichlet energy  $\mcI_1=|\nabla Q|^2=Q_{ij,k}Q_{ij,k}$.  The   
bulk energy density is necessarily  a function of the spectral invariants $\tr (Q^2)$ and $\tr (Q^3)$ (recall that $\textrm{tr}\, Q=0$), and is typically taken to be (see, eg,  \cite{mottram2014introduction})
\be\label{def:fb}
f_B(Q) = -\frac{a^2}{2}\tr(Q^2)-\frac{b^2}{3}\tr(Q^3)+\frac{c^2}{4}\left(\tr(Q^2) \right)^2.
\ee
Taking the elastic energy to be the Dirichlet energy, in either the Oseen-Frank and Landau-de Gennes theories, leads to considerable simplifications.  The second-order differential operator in the Euler-Lagrange equation is just the scalar Laplacian, which does not couple the components of $n$ and $Q$ respectively  (coupling occurs only through first-order 
or zeroth-order 
derivative terms respectively).  
This enables techniques from  scalar PDEs to be carried over to the system of Euler-Lagrange equations.  

However, there is no physical justification for imposing $O(3)$-invariance separately in physical space and order-parameter space; only joint invariance is required.  For the Oseen-Frank theory, requiring  $\sigma(n,\nabla n)$ to be invariant under $n(x) \mapsto Pn(P^{-1}x)$ implies that the energy density is of the following form  (see, eg,  Section $3.2$ of Chapter $3$ in \cite{Virga}) :
\be\label{sigmaOF}
\sigma(n,\nabla n)=K_1|\nabla\cdot n|^2+K_2 |n\cdot(\nabla\times n)|^2+K_3|n\times (\nabla\times n)|^2+(K_2+K_4)\left[\textrm{tr}(\nabla n)^2-(\nabla\cdot n)^2\right].
\ee
 The Oseen-Frank energy \eqref{sigmaOF} is found to explain a wide variety of equilibrium phenomena in nematic liquid crystals. (Requiring only $SO(3)$-invariance allows for an additional term linear in $\nabla n$, which can describe cholesteric liquid crystals).
The {\it elastic constants} $K_1,K_2,K_3$ and $K_4$, which are free parameters in the energy density, have been measured in experiments for a variety of  materials as well as calculated through numerical simulations.    

The  term $\textrm{tr}(\nabla n)^2-(\nabla\cdot n)^2$ in \eqref{sigmaOF}, called the saddle-splay term, plays the role of a null Lagrangian; it can be expressed as the 
pure divergence $\nabla\cdot((n\cdot \nabla) n - (\nabla \cdot n) n)$, whose volume integral depends only on  $n$ and its tangential derivatives on the boundary.  
If the boundary conditions are fixed,  the saddle-splay term can be ignored, at least for the purposes of energy minimisation.
The Dirichlet energy is recovered by taking $K_1 = K_2 =K_3$ and  $K_4=0$. 
In general, the elastic constants should be chosen so that the energy density is bounded from below.  The resulting constraints, first derived by Ericksen\cite{ericksen_inequalities}, are given by
\be \label{eq: Ericksen ineq}
2K_1 > K_2 + K_4 > 0.
\ee

Similarly for the Landau-de Gennes theory, $O(3)$-invariance requires only that the energy density be invariant under $Q(x) \mapsto PQ(P^{-1} x) P^t$.  The conditions on the bulk energy remain the same, but the elastic energy $f_E(Q,\nabla Q)$ may assume a much more general form. 
Using methods of representation theory and invariant theory, one can construct the linearly independent,  $O(3)$-invariant polynomials quadratic in $\nabla Q$ and of order $m\ge 0$ in $Q$. 
(If one requires only $SO(3)$-invariance, there are   additional invariants linear in $\nabla Q$, which describe cholesteric liquid crystals.)
The simplest  of the invariants quadratic in $\nabla Q$ are those which are independent of $Q$ (ie, $m=0$).  
There are just three, namely the Dirichlet energy $\mcI_1$ and the following two:
$$\mcI_{2}= Q_{ik,j}Q_{ij,k},\quad\mcI_3=Q_{ij,j} Q_{ik,k}.$$
The elastic energy assumes the form
\be \label{eq: f_E}
f_E(\nabla Q) = \frac{L_1}{2}Q_{ij,k}Q_{ij,k} +\frac{L_2}{2} Q_{ik,j}Q_{ij,k}+\frac{L_3}{2}Q_{ij,j}Q_{ik,k} 
\ee
with  elastic constants $L_1$, $L_2$ and $L_3$ as free parameters.

It is well known that $\mcI_2$ and $\mcI_3$  can be combined to give a null Lagrangian, analogous to the saddle-splay term in the Oseen-Frank energy.  As this fact appears to be difficult to find in the literature, we state the result formally here and provide a proof in the Appendix.
\begin{lemma} {\bf [Null Lagrangian]} 
\label{lemma:nulllagrangian} Let $\Omega\subset\RR^d$, $d=2,3$, be a domain with $C^1$ boundary. Let $Q\in H^1(\Omega)$ and furthermore, if $d=2$, assume that $\frac{\partial Q}{\partial x_3}\equiv 0$. Then: 
$$\int_\Omega ( \mcI_2-\mcI_3) \,dx=\int_\Omega Q_{ik,j} Q_{ij,k} \, dx -\int_\Omega Q_{ij,j}Q_{ik,k}\,dx$$ depends only on $Q|_{\partial\Omega}\in H^{\frac{1}{2}}(\partial\Omega)$.
\end{lemma}
In general,
the elastic constants $L_1$, $L_2$ and $L_3$ should be chosen so that the energy density is bounded from below.
The resulting constraints  are obtained, for instance, in \cite{Longa-several}, using invariant theory arguments. A different approach allows to obtain the same result by standard methods. We have thus the following lemma, whose proof is given in the Appendix.
\begin{lemma}{\bf [Coercivity conditions]}
\label{lemma:coercivity} Let $\Omega\subset\mathbb{R}^d$, $d=2,3$. Let  $Q\in H^1(\Omega;\mcS_0)$, and furthermore, if $d=2$, assume $\frac{\partial Q}{\partial x_3}\equiv 0$.
\begin{enumerate}

\item 
There exists $\mu_0>0$, independent of $Q$, such that for almost all $x\in \Omega$

\be
\frac{L_1}{2}Q_{ij,k}Q_{ij,k} +\frac{L_2}{2} Q_{ik,j}Q_{ij,k}+\frac{L_3}{2}Q_{ij,j}Q_{ik,k}\ge \mu_0|\nabla Q|^2
\ee if and only if
\be \label{sufcc}
L_1+L_2>0, \ 2L_1-L_2>0, \  L_1+\frac{L_2}{6}+\frac{5}{3}L_3>0.
\ee

\item  Let $M:=\frac{L_1+L_2}{2}$ and $Q_b\in H^\frac{1}{2}(\partial\Omega;\mcS_0)$. Assume that

\be \label{eq: coercivity inequalities}
L_1 + \frac43 M>0,\ \  L_1 >0.
\ee 

Then there exists $\mu_1>0$ so that for any $Q\in H^1(\Omega;\mcS_0)$  with $Q=Q_b$ we have:

\be
\int_\Omega \frac{L_1}{2}Q_{ij,k}Q_{ij,k} +MQ_{ij,j}Q_{ik,k}\ge \int_\Omega \mu_1 |\nabla Q|^2\,dx+\mathcal{B}
\ee with $\mathcal{B}$ a constant depending only on $Q_b$, $L_2$ and $L_3$, but independent of $Q$.

\end{enumerate}
\end{lemma}

One can  connect the $Q$-tensor theory and the Oseen-Frank theory by taking $Q$ to be of the uniaxial form
 \be\label{Qlimit}
 Q=s_+(n\otimes n-\frac{1}{3}I),
 \ee 
 with $n \in S^2$ and fixed scalar order parameter $s_+$.  Some rigorous justifications for this are provided in \cite{Ma-Za, Ng-Za} with a physical perspective available in \cite{gartland2015scalings}.  The constraints are motivated by the fact that uniaxial $Q$-tensors with $s_+$  appropriately chosen are precisely the (degenerate) global minimizers of the bulk potential $f_B$ given by \eqref{def:fb}. 
In a regime where, after suitable nondimensionalization,  the elastic constants may be regarded as small, minimizers of the Landau-de Gennes energy approach the uniaxial form  \eqref{Qlimit} nearly  everywhere.
 
For uniaxial $Q$-tensors \eqref{Qlimit}, the transformation $Q(x) \mapsto P Q(P^{-1} x) P^t$ is equivalent to $n(x) \mapsto Pn(P^{-1} x)$ for $P \in O(3)$.  It follows that the Landau-de Gennes elastic energy density $f_E$ evaluated on uniaxial $Q$-tensors  yields an expression in $n$ and $\nabla n$ of the Oseen-Frank form \eqref{sigmaOF}.
Straightforward calculations lead to the following identifications of the resulting Oseen-Frank elastic constants in terms of the Landau-de Gennes elastic constants (for simplicity $s_+=1$):
\be \label{eq: OF and LDG} K_1 = K_3 = L_1 + \frac12(L_2 + L_3), \quad K_2 = L_1, \quad K_2 + K_4 = L_1 + \frac12 L_2. \ee
Thus, for uniaxial $Q$ tensors with fixed scalar order parameter, the Landau-de Gennes elastic energy \eqref{eq: f_E} reduces to a restricted form of the Oseen-Frank energy, in which $K_1 = K_3$ (that there must be a restriction is already evident from the fact that the Oseen-Frank elastic energy \eqref{sigmaOF} has four elastic constants while the Landau-de Gennes elastic energy \eqref{eq: f_E} has only three).
The restriction can be overcome by allowing the next-most-simple set of terms in the Landau-de Gennes elastic energy, namely terms quadratic in $\nabla Q$ and linear in $Q$ (ie, $m = 1$ above).
This introduces additional complications.  First, there are six such terms (including an additional null Lagrangian), leading to a surfeit of elastic constants as compared to the Oseen-Frank theory  rather than a deficit. 
Second, the resulting elastic energy density is unbounded below.   The latter problem can be addressed by introducing yet more elastic terms, say quadratic in both $Q$ and $\nabla Q$, and/or by modifying the bulk energy \eqref{def:fb} so as to be infinite outside a compact subset of $Q$-tensors (such modifications are independently motivated by considerations of the microscopic model -- see \cite{Ball-Maj}, \cite{singpotphys} in $d=3$ and \cite{onsager-spatial_variation} in $d=2$). 
In any case, we may conclude that if there are phenomena well described by the Oseen-Frank theory in which the fidelity of the description  (qualitative or quantitative)  depends  on having $K_1 \neq K_3$,  then the several-elastic-constant Landau-de Gennes elastic energy \eqref{eq: f_E}, which is the subject of this paper,  may not provide an adequate description.

 \bigskip
In the remainder of this paper we will focus on studying the effects induced by the presence of the anisotropic elastic energy terms, $\mcI_2$ and $\mcI_3$, on a basic model problem that is well understood in the case of one elastic constant, namely the two-dimensional  k-radially symmetric solutions, studied for instance in \cite{DRSZ,INSZ2d1,INSZ2d2}.

\subsection{Mathematical formulation of the problem and main results}
 As we will be imposing Dirichlet boundary conditions, in view of Lemma~\ref{lemma:coercivity} it suffices to restrict the elastic energy to have just two independent constants, namely  $L := L_1$ and $M := (L_2+L_3)/2$. We therefore consider the following form of the Landau-de Gennes energy:
\begin{align}\label{LDG}
\mcF [Q; \Omega]= \int_{\Omega} \Big[\frac{L}{2}Q_{ij,k}Q_{ij,k} +M \, Q_{ij,k}Q_{ik,j}+ f_B(Q) \Big]\,dx,
\end{align} 
where 
$L$, $M$ 
are the independent elastic constants with $L>0$, $L+\frac43 M > 0$,
and $Q \in H^1(\Omega, \mcS_0)$ takes values in the space of $Q$-tensors,
\be\label{def:s0}
\mcS_0\defeq\{ Q\in \RR^{3\times 3},\, Q=Q^t,\,\tr (Q)=0\},
\ee
The  bulk potential $f_B(Q)$ is  taken to be of the form \eqref{def:fb} where $a^2, b^2 \geq 0$ and $c^2 > 0$ are material constants. 
Here and in what follows, summations are taken over $i,j,k=1,2,3$.   Since we work in a two-dimensional domain, we assume throughout $Q_{ij,3}\equiv 0, \forall i,j=1,2,3$. The relationship between the $2D$ and $3D$ solutions is discussed in the Appendix, Section~\ref{sec:2d3d}.

We are interested in studying  two-dimensional point defects in liquid crystals and their characteristic symmetry features. The simplest and most generic liquid crystal point defects are obtained as critical points of \eqref{LDG} in the domain $\Omega= B_R \subset\RR^2$ (the ball  of radius $R>0$ centered at the origin) under the following boundary conditions  (for more details see \cite{DRSZ}):
\be\label{BC1}
Q(x)=Q_k(x) \equiv s_+ \left( n (x) \otimes n (x) -\frac{1}{3} I \right), \ x\in\partial B_R(0),
\ee
where $x = (R\cos\varphi,  R\sin \varphi)$, 
\be\label{nn}
n(x) =\left(\cos ({\textstyle\frac{k}{2}} \varphi) , \sin ({\textstyle\frac{k}{2}} \varphi) , 0\right), \ k \in \ZZ \setminus \{0\}, 
\ee
and 
\be\label{s+}
s_+ = \frac{ b^2 + \sqrt{b^4+24 a^2 c^2}}{4 c^2}.
\ee The value of $s_+$ is chosen  so that $Q_k$ minimizes $f_B(Q)$. 

\vskip 0.1cm

The critical points of the energy \eqref{LDG} satisfy the following Euler-Lagrange equations:
\bea \label{eq:EL}
L \Delta Q_{ij} &+ M\left( \partial_j \partial_k Q_{ik} + \partial_i \partial_k Q_{jk}  - \frac{2}{3} \partial_l \partial_k Q_{lk} \delta_{ij} \right) \non \\  &= -a^2 Q_{ij} -b^2\left( Q_{ik}Q_{kj} - \frac{1}{3} |Q|^2 \delta_{ij} \right) + c^2 Q_{ij} |Q|^2 \, , i,j =1, 2, 3 \ \ \hbox{ in } B_R,
\eea 
subject to the boundary conditions \eqref{BC1}. 
We note that terms $-\frac{2M}{3}\partial_l \partial_k Q_{lk} \delta_{ij}$ and $\frac{b^2}{3}|Q|^2 \delta_{ij}$ account for the constraint $\tr(Q)=0$. 

\bigskip
We are interested in studying the critical points compatible with the symmetry of the problem \eqref{eq:EL}, \eqref{BC1}. 

\begin{definition}
\label{def:fullyk_rad}
For $k \in \ZZ \setminus \{ 0 \}$, we say that a Lebesgue measurable map $Q: B_R \to \mcS_0$ is {\bf generally} $k$-{\bf radially symmetric} if the following condition holds for almost every $x=(x_1, x_2)\in B_R$:
\vskip 0.2cm
$$
Q\bigg(P_2\big({\mathcal R}_2 (\psi)  \tilde x\big)\bigg)= {\mathcal R}_k (\psi) Q(x) {\mathcal R}_k^t (\psi) , \ \textrm{for almost  every } \psi \in \RR, 
$$ where $\tilde x=(x_1,x_2,0)$, $P_2:\RR^3\to \RR^2$ is the projection given by $P_2(x_1,x_2,x_3)=(x_1,x_2)$, and
\begin{equation}
\label{eq: R_k }
{\mathcal R}_k (\psi) := \left(\begin{array}{ccc}\cos(\frac{k}{2}\psi) & -\sin(\frac{k}{2}\psi) & 0 \\\sin(\frac{k}{2}\psi) & \cos(\frac{k}{2}\psi) & 0 \\0 & 0 & 1\end{array}\right)
\end{equation} 
represents the rotation about  $e_3$ by angle $\frac{k}{2}\psi$.
\end{definition}
 We recall that in \cite{INSZ2d1} we studied a class of critical points of Landau-de Gennes energy with restricted symmetry, namely:
\begin{definition}
\label{def:k_rad}
For $k \in \ZZ \setminus \{ 0 \}$, we say that a Lebesgue measurable map $Q: B_R \to \mcS_0$ is  $k$-{\bf radially symmetric} if it is generally $k$-radially symmetric and the following additional condition holds:
\begin{itemize}
\item [{\bf(H1)}] The vector $e_3 =(0,0,1)$ is an eigenvector of $Q(x)$, for almost all $x\in B_R$.
\end{itemize}
\end{definition}
\begin{remark}
If $k$ is an odd integer, then a generally $k$-radially symmetric  $Q \in H^1 (\Omega, \mcS_0)$  automatically verifies {\bf (H1)},  so there is no difference between the two types of symmetries (see Proposition $2.1$ in \cite{INSZ2d1}). 
\end{remark}

In the case 
$M =0$, {\bf all} $k$-radially symmetric critical points of Landau-de Gennes energy \eqref{LDG} with boundary condition \eqref{BC1} have the following form \cite{DRSZ, INSZ2d1, INSZ2d2}
\be\label{anY}
Y_k= u(r) \sqrt{2}\left(n\otimes n-\frac{1}{2}I_2\right) + v(r) \sqrt{\frac{3}{2}}\left(e_3\otimes e_3-\frac{1}{3}I_3\right),
\ee
where $e_i$, $i=1,2,3$ are the standard cartesian basis vectors, $I_2 =e_1 \otimes e_1 + e_2 \otimes e_2$, $I_3=I_2+e_3\otimes e_3$  and $(u,v)$ satisfies the following system of ODEs:
\begin{align}\label{ODEsystem}
u''+\frac{u'}{r}-\frac{k^2u}{r^2} &=\frac{u}
{L}\left[-a^2+\sqrt{\frac{2}{3}} b^2 v+c^2\left( u^2+ v^2\right)\right],\nonumber\\
v''+\frac{v'}{r}&=\frac{v}
{L}\left[-a^2-\frac{1}{\sqrt{6}}b^2 v+c^2\left( u^2+ v^2\right) \right] + \frac{1}{\sqrt{6} L} b^2 u^2
\end{align}
with boundary conditions 
\be\label{bdrycond}
u(0)=0,\ v'(0)=0, \ u(R)=\frac{1}{\sqrt{2}} s_+,\,\,\,v(R)=-\frac{1}{\sqrt{6}}s_+.
\ee

Moreover, it has been shown that
\begin{itemize}
\item in the case $b=0$, the $k$-radially symmetric critical points on $B_R$ ($R< \infty)$ are global minimizers of \eqref{LDG}, \eqref{BC1} for any $k \in \ZZ \setminus \{0\}$ (see \cite{DRSZ});
\item in the case $b \neq 0$ all $k$-radially symmetric critical points on $\R^2$ are unstable for $k\neq \pm 1$ (see \cite{INSZ2d1}),  while for $k = \pm 1$, there exists solutions  \eqref{anY} with $u>0$ and $v<0$ on $(0,R)$ with $R \leq \infty$ that are stable.
\end{itemize}

The main aim of this paper is to study the existence and  behaviour  of  {\it generally $k$-radially symmetric} critical points  as well as nonsymmetric local minimisers of the Landau-de Gennes energy \eqref{LDG} with boundary conditions \eqref{BC1} for general $L$ and $M$ subject to the coercivity 
conditions~\eqref{eq: coercivity inequalities}.  We also investigate the symmetry-breaking of radial critical points for small nonzero $M$.

\vskip 0.2cm

We denote, for $\varphi \in [0,2\pi)$,
\[
n=\left(\cos ({\textstyle\frac{k}{2}} \varphi) , \sin ({\textstyle\frac{k}{2}} \varphi) , 0\right), \,
m=\left(-\sin({\textstyle\frac{k}{2}} \varphi),\cos({\textstyle\frac{k}{2}} \varphi),0\right) .
\]
We endow the space $\mcS_0$ of $Q$-tensors with the scalar product $$Q\cdot \tilde Q=\tr(Q\tilde Q)$$ and 
for any $\f\in [0, 2\pi)$, we define the following orthonormal basis in $\mcS_0$:
\begin{gather}
E_0=\sqrt{\frac{3}{2}}\left(e_3\otimes e_3-\frac{1}{3}I_3\right),\\
E_1 =E_1(\f)=\sqrt{2}\left(n\otimes n-\frac{1}{2}I_2\right),\quad E_2=E_2(\f)=\frac{1}{\sqrt{2}}\left(n\otimes m+m\otimes n\right), \non \\
E_3 =\frac{1}{\sqrt{2}}(n\otimes e_3+e_3\otimes n),\quad E_4=\frac{1}{\sqrt{2}}\left(m\otimes e_3+e_3\otimes m\right).
\end{gather}
It is straightforward to check that
\be
\label{deriv_E12}
\frac{\partial E_1}{\partial \f}=kE_2, \quad  \frac{\partial E_2}{\partial \f}=-kE_1, \quad 
\frac{\partial E_3}{\partial \f}=\frac{k}{2} E_4, \quad  \frac{\partial E_4}{\partial \f}=-\frac{k}{2}E_3.
\ee

\begin{remark}
The above basis is relevant only for the case when $k$ is even.  For $k$ odd, $E_3$ and $E_4$ are not $2\pi$-periodic in $\varphi$.  In this case,
a good basis can be obtained as in \cite{INSZ2d1} by replacing $E_3$ with $\tilde E_3=\frac{1}{\sqrt{2}}(e_1\otimes e_3+e_3\otimes e_1)$ and $E_4$ with $\tilde E_4=\frac{1}{\sqrt{2}}(e_2\otimes e_3+e_3\otimes e_2)$. 
\end{remark}

Our main analytical result provides conditions for the existence of generally $k$-radially symmetric critical points of the Landau-de Gennes energy.
\begin{theorem}
\label{thm:radialsym}
Assume that $k \in \ZZ \setminus\{0\}$, 
$L>0$, $M \neq 0$, and $L + \frac43 M >0$. 
\begin{enumerate}
\item If $k\neq2$ then there exist no generally $k$-radially symmetric critical points of the Landau-de Gennes energy \eqref{LDG} subject to boundary conditions \eqref{BC1}.
\item If $k=2$ then  $Q \in H^1(B_R, \mcS_0)$ is a generally $k$-radially symmetric critical point of the Landau-de Gennes energy \eqref{LDG} satisfying the boundary conditions \eqref{BC1} if and only if
\be\label{5comp}
Q(x) = \sum_{i=0}^4 w_i(r) E_i,
\ee
where $w_i \in C^\infty(0,R)$ satisfy the following system of ODE's:
\bea\label{sys5}
(L&+ M/3 ) \left(w_0''+\frac{w_0'}{r}\right) - \frac{M}{\sqrt{3}} \left( w_1'' + \frac{3 w_1'}{r} \right)  \non \\ 
 &=   w_0 \left( -a^2 - \frac{ b^2}{\sqrt{6}} w_0 + c^2 \left(\sum_{i=0}^4 w_i^2\right) \right) + \frac{ b^2}{\sqrt{6}} (w_1^2+w_2^2) - \frac{ b^2}{2\sqrt{6}} (w_3^2+w_4^2), \non \\ 
 (L&+M) \left(w_1''+\frac{w_1'}{r}-\frac{4w_1}{r^2}\right) - \frac{M}{\sqrt{3}} \left( w_0'' - \frac{w_0'}{r} \right) \non \\ 
 &= w_1 \left( -a^2 + \frac{2 b^2}{\sqrt{6}} w_0 + c^2 \left(\sum_{i=0}^4 w_i^2\right) \right) - \frac{ b^2}{2\sqrt{2}} (w_3^2-w_4^2),\non \\
 (L&+M) \left( w_2'' + \frac{w_2'}{r} - \frac{4 w_2}{r^2} \right)  =  w_2  \left( -a^2 + \frac{2 b^2}{\sqrt{6}} w_0 + c^2 \left(\sum_{i=0}^4 w_i^2\right) \right) - \frac{ b^2}{\sqrt{2}} w_3 w_4 ,\\
(L&+M) \left( w_3'' + \frac{w_3'}{r} - \frac{w_3}{r^2} \right)  =  w_3  \left( -a^2 - \frac{b^2}{\sqrt{6}} w_0 - \frac{b^2}{\sqrt{2}} w_1+ c^2 \left(\sum_{i=0}^4 w_i^2\right)  \right) - \frac{ b^2}{\sqrt{2}} w_2 w_4, \non \\
L& \left( w_4'' + \frac{w_4'}{r} - \frac{w_4}{r^2} \right)  =  w_4  \left( -a^2 - \frac{b^2}{\sqrt{6}} w_0 + \frac{b^2}{\sqrt{2}} w_1+ c^2 \left(\sum_{i=0}^4 w_i^2\right)  \right) - \frac{ b^2}{\sqrt{2}} w_2 w_3, \non 
\eea
subject to the boundary conditions 
\bea\label{BCode5}
&  
w_0'(0)=0, w_1(0)=0, \ w_2(0)=0, \ w_3(0) =0, \ w_4(0)=0, \non \\ 
&w_0(R) = -\frac{s_+}{\sqrt{6}},  \ w_1(R)=\frac{s_+}{\sqrt{2}},\  w_2(R)=0, \ w_3(R)=0, \ w_4(R)=0.
\eea
\end{enumerate}
\end{theorem}

Using this result,  in Section~\ref{sec3} we carry out a numerical  study of the generally $k$-radially symmetric  solutions for $k=2$.  Depending on 
the elastic constant $M$, the material parameter $b$ and the radius $R$, we  observe three types of  solutions, which  we  classify in terms of the number of non-vanishing coefficients  $w_i$, as follows:
\begin{itemize}
\item {\it two-component solutions}, where $Q = w_0(r) E_0 + w_1(r) E_1$;
\item {\it  three-component solutions}, where $Q = w_0(r) E_0 + w_1(r) E_1 + w_3(r) E_3$;
\item {\it  five-component solutions} where $Q =\sum_{i=0}^4 w_i(r) E_i$.
\end{itemize}
These solutions appear as minimizers (local or global) of the full Landau-de Gennes energy in certain parameter regimes. We determine the profiles and stability properties of these  solutions, and discuss a number of their  characteristic features, including symmetries, the monotonicity and signs of the radial components.  
We also find  non-radially-symmetric solutions with $\mathbb{Z}_2$ symmetry,  qualitatively similar to  previously reported solutions \cite{baumanphilips, disclinationsnumerics}, exhibiting two separated index-one-half defects.  
By comparing the (numerically determined) energies of all these solutions, we obtain the global phase diagram.

In the Section~\ref{sec4} we investigate how general $k$-radial symmetry is broken  for small $M \neq 0$ in the physically interesting cases $k=\pm1$. Using formal asymptotics, we characterise the leading-order symmetry-breaking terms as solutions of a system of linear inhomogeneous ODE's generated by the second variation of the Landau-de Gennes energy.  These asymptotic results are supported by numerical calculations.

In the last section of the paper we present a number of open problems (both analytical and numerical)  motivated by our numerical investigations.  It is hoped that these will  stimulate the development of new  techniques for  addressing the physically relevant but mathematically challenging Landau-de Gennes model with several elastic constants. Finally, in the Appendix we cover a number of technical questions necessary for the paper.


\section{Radially symmetric critical points}

This section is devoted to an analytical study of generally $k$-radially symmetric critical points of the Landau-de Gennes energy \eqref{LDG} subject to the boundary conditions \eqref{BC1}. In particular, we show that if $M \neq 0$, then generally $k$-radially symmetric critical points do not exist if $k \neq 2$. In the case $k=2$ we derive the system of ODEs satisfied by generally $k$-radially symmetric critical points and using variational methods prove the existence of a solution for this system. Finally, we study the interesting limit $M \to \infty$.

We begin with the proof of  Theorem~\ref{thm:radialsym} about existence of generally $k$-radially symmetric critical points.

\bigskip
\noindent {\bf Proof of Theorem~\ref{thm:radialsym}.}  We know that any generally $k$-radially symmetric $Q$-tensor $Q \in H^1(B_R; \mcS_0)$ can be represented in the  form 
$$
Q= 
 \sum_{i = 0}^4
w_i(r) E_i,
$$ 
where $w_0 \in H^1((0,R); r\,dr)$ and $w_i  \in H^1((0,R); r\,dr)\cap  L^2\left((0,R); \frac{1}{r}\,dr\right)$ for $1 \leq i \leq 4$. Moreover, for $k$ odd, the components $w_3,w_4$ vanish, see Proposition $2.1$ in  \cite{INSZ2d1}. 

Similarly as in \cite{INSZ2d1} we prove now that $w_i(0) = 0$ for $1\le i\le 4$. Since  $w_i\in H^1((0,R); r\,dr)\cap  L^2\left((0,R); \frac{1}{r}\,dr\right)$ we have that that $w_i$ is continuous on $(0,R)$ and for $r_1, r_2\in (0,R)$:
$$|w_i^2(r_2)-w_1^2(r_1)|=2\bigg|\int_{r_1}^{r_2}w_iw_i'\, dr\bigg|\leq 2\bigg(\int_{r_1}^{r_2}{w_i^2}\, \frac{dr}{r}\bigg)^{1/2} 
\bigg(\int_{r_1}^{r_2} (w'_i)^2\, rdr\bigg)^{1/2}.$$ Since the right hand side converges to zero as $|r_2-r_1|\to 0$, it follows that $w_i$ is continuous up to $r=0$. Combined again with $w_i\in L^2\left((0,R); \frac{1}{r}\,dr\right)$, we conclude that $w_i(0)=0$. 

We would like to understand under what conditions the generally $k$-radially symmetric $Q$-tensor satisfies the full Euler-Lagrange equations \eqref{eq:EL} with boundary conditions \eqref{BC1}. The nonlinear part, contained in the right-hand side of the Euler-Lagrange equations \eqref{eq:EL}, becomes
\bea\label{non5}
-a^2 Q - b^2& \left(Q^2 -\frac{1}{3} |Q|^2 I \right) + c^2 Q |Q|^2 =  \\ 
& \left[ w_0 \left( -a^2 - \frac{ 2 b^2}{\sqrt{6}} w_0 + c^2 \left(\sum_{i=0}^4 w_i^2\right) \right) + \frac{ b^2}{\sqrt{6}}  \left(\sum_{i=0}^4 w_i^2\right)  - \frac{ 3 b^2}{2\sqrt{6}} (w_3^2+w_4^2) \right] E_0 \non \\ 
&+ \left[  w_1  \left( -a^2 + \frac{2 b^2}{\sqrt{6}} w_0 + c^2 \left(\sum_{i=0}^4 w_i^2\right) \right) - \frac{ b^2}{2\sqrt{2}} (w_3^2-w_4^2) \right] E_1 \non \\ 
&+ \left[ w_2  \left( -a^2 + \frac{2 b^2}{\sqrt{6}} w_0 + c^2 \left(\sum_{i=0}^4 w_i^2\right) \right) - \frac{ b^2}{\sqrt{2}} w_3 w_4 \right] E_2 \non \\
&+ \left[ w_3  \left( -a^2 - \frac{b^2}{\sqrt{6}} w_0 - \frac{b^2}{\sqrt{2}} w_1+ c^2 \left(\sum_{i=0}^4 w_i^2\right)  \right) - \frac{ b^2}{\sqrt{2}} w_2 w_4 \right] E_3 \non \\
&+ \left[  w_4  \left( -a^2 - \frac{b^2}{\sqrt{6}} w_0 + \frac{b^2}{\sqrt{2}} w_1+ c^2 \left(\sum_{i=0}^4 w_i^2\right)  \right) - \frac{ b^2}{\sqrt{2}} w_2 w_3 \right] E_4, \non
\eea
with $w_3=w_4=0$ in the case of odd $k$.

The elastic part, contained in the left-hand side of the Euler-Lagrange equations \eqref{eq:EL},   provides an elliptic operator under the assumed conditions on $L$ and $M$. Standard arguments (see for instance \cite{Tim-Gart}) show that the solutions are smooth. Then  one can easily check that $$w_i \in C^\infty(0,R)$$ for $i=0,\dots 4$ as claimed.

 Furthermore,  since $Q$ is smooth in $B_R$, we obtain that $\textrm{tr}(Q E_0)$ is smooth in $B_R$ and in particular along the line $\{(x,0); -R\le x\le R\}$. Noting that $w_0(r)=\textrm{tr}(Q(r,0)E_0)$ for  $r\geq0$ and taking into account that $Q$ is radially symmetric, we can trivially extend $w_0$ to $(-R,R)$ in an even and smooth way. Therefore $w_0 \in C^\infty([0,R))$ and $w_0'(0) = 0$.

Furthermore the elastic part generates a linear operator. Therefore, we can study the action of the elastic part on the components $w_i(r) E_i$ 
separately. For convenience, 
we define the  operator 
\be\label{Lop}
 {\mathcal L} Q_{ij}= \partial_j \partial_k Q_{ik} + \partial_i \partial_k Q_{jk}  - \frac{2}{3} \partial_l \partial_k Q_{lk} \delta_{ij} .
\ee
We notice that
\be\label{el012}
(L \Delta + M {\mathcal L} ) \left(\sum_{i=0}^2 w_i(r) E_i \right) \in \hbox{\rm span} \{ E_0, E_1, E_2\}
\ee
and
\be\label{el34}
(L \Delta + M {\mathcal L} ) \left(\sum_{i=3}^4 w_i(r) E_i \right) \in \hbox{\rm span} \{ E_3, E_4\},
\ee
with the convention that $w_3=w_4=0$ when $k$ is odd.

Let us first compute the elastic part for \eqref{el012}. We note that for any $k \in \ZZ \setminus \{0\}$ (see \cite{DRSZ})
\bea\label{deltatopolar}
\Delta \left(\sum_{i=0}^2 w_i(r) E_i \right) =\left(w_0''+\frac{w_0'}{r}\right)E_0 + \left(w_1''+\frac{w_1'}{r}-\frac{k^2w_1}{r^2}\right)E_1+ \left(w_2''+\frac{w_2'}{r}-\frac{k^2w_2}{r^2}\right)E_2.
\eea
After a lengthy but straightforward calculation (see Appendix, Section~\ref{sec:symmbreaking}) we also obtain 
\bea\label{Lop3}
 {\mathcal L} (w_0(r) E_{0} ) &= \frac{1}{3} \left( w_0'' + \frac{w_0'}{r} \right) E_0 - \frac{1}{\sqrt{3}} \left( w_0'' - \frac{w_0'}{r} \right) \left( E_1 \cos((k-2) \varphi) - E_2 \sin((k-2) \varphi) \right), \non \\
{\mathcal L} (w_1(r) E_{1} ) &= \left( w_1'' + \frac{w_1'}{r} - \frac{k^2 w_1}{r^2} \right) E_1 - \frac{E_0}{\sqrt{3}} \left( w_1'' + \frac{(2k-1) w_1'}{r} +\frac{k (k-2) w_1}{r^2} \right) \cos ((k-2)\varphi),  \\
{\mathcal L} (w_2(r) E_{2} ) &= \left( w_2'' + \frac{w_2'}{r} - \frac{k^2 w_2}{r^2} \right) E_2 + \frac{E_0}{\sqrt{3}} \left( w_2'' + \frac{(2k-1) w_2'}{r}+ \frac{k (k-2) w_2}{r^2} \right) \sin ((k-2)\varphi). \non
\eea
Since $E_0$, $E_1$ and $E_2$ are 
linearly independent, we 
see that the only possibility to obtain a closed system of  ordinary differential equations for $w_0$, $w_1$ and $w_2$ is to take $k=2$. Therefore, even without computing 
the terms in \eqref{el34}, we may conclude that for any $k \neq 2$, $M \neq 0$  there is no generally $k$-radially symmetric solution of \eqref{eq:EL}. 

\bigskip
Now we want to investigate the case $k=2$. We already computed the term $(L \Delta + M {\mathcal L}) \left( \sum_{i=0}^2 w_i(r) E_i \right)$ and we are left with finding $(L \Delta + M {\mathcal L}) \left( \sum_{i=3}^4 w_i(r) E_i \right)$. After a straightforward calculation we obtain
\be\label{Lop5}
(L \Delta + M {\mathcal L}) \left( \sum_{i=3}^4 w_i(r) E_i \right) = (L+M) E_3\left(w_3''+\frac{w_3'}{r}-\frac{w_3}{r^2}\right) +   L E_4 \left(w_4''+\frac{w_4'}{r}-\frac{w_4}{r^2}\right) .
\ee
Combining \eqref{non5}, \eqref{deltatopolar}, \eqref{Lop3}, \eqref{Lop5} we observe that the Euler-Lagrange equations become \eqref{sys5}. Using boundary conditions \eqref{BC1} for the PDE system we obtain the  boundary conditions for the ODE at $R$.  

In order to prove existence of the solution for the problem \eqref{sys5}, \eqref{BCode5} we minimize the corresponding energy functional
\bea\label{def:mcR5}
\mcE(\{w_i\})&= \int_0^R \left\{ \frac{L}{2} \left(|w_1'|^2 + |w_0'|^2 + \frac{4}{r^2} w_1^2 \right) + \frac{M}{6}  \left| \sqrt{3} w_1' -w_0' + \frac{2\sqrt{3}}{r} w_1 \right|^2 \right\} r \, dr \non \\ 
 &+ \int_0^R \left\{ \frac{L+M}{2} \left(|w_2'|^2  + \frac{4}{r^2} w_2^2 + |w_3'|^2  + \frac{1}{r^2} w_3^2 \right) + \frac{L}{2} \left( |w_4'|^2  + \frac{1}{r^2} w_4^2\right)  \right\} r \, dr \non \\
 &+ \int_0^R \left\{ \left( -\frac{a^2}{2}+\frac{c^2}{4} \left(\sum_{i=0}^4 w_i^2\right) \right) \left( \sum_{i=0}^4 w_i^2 \right) \right.\\
 &\left. - \frac{b^2\sqrt{6}}{36}\left( 2w_0^3-6w_0(w_1^2+w_2^2) +3w_0 (w_3^2+w_4^2) + 3\sqrt{3} w_1(w_3^2-w_4^2) + 6 \sqrt{3} w_2 w_3 w_4 \right) \right\} r \, dr \non
\eea
defined on the admissible set $S \times (S_w)^3$ where $S$ and $S_w$ are defined as
\be\label{SR}
S = \left\{ (w_0,w_1) :  [0,R] \to \RR^2 \, \Big | \, \sqrt{r} w_0', \sqrt{r} w_1',  \sqrt{r} w_0, \frac{w_1}{\sqrt{r}} \in L^2(0,R), \, w_0(R)= -\frac{s_+}{\sqrt{6}} , w_1(R)=\frac{s_+}{\sqrt{2}}\right\}
\ee
and 
\be\label{SRw}
S_w = \left\{ w :  [0,R] \to \RR\, \Big | \, \sqrt{r} w', \frac{w}{\sqrt{r}} \in L^2(0,R), \, w(R)=0 \right\}.
\ee
It is clear that for $M\geq0$ the above energy is coercive and therefore existence of a minimizer follows by standard arguments. In the case $M<0$ we can rewrite the energy as 
\bea\label{def:mcR5n}
\mcE(\{w_i\})&= \frac{2 M s_+^2}{3} + \int_0^R \left\{ \left(\frac{L}{2} + \frac{4M}{6} \right) \left(|w_1'|^2 + |w_0'|^2 + \frac{4}{r^2} w_1^2 \right) - \frac{M}{6}  \left| \sqrt{3} w_0' +w_1' + \frac{2}{r} w_1 \right|^2 \right\} r \, dr \non \\ 
 &+ \int_0^R \left\{ \frac{L+M}{2} \left(|w_2'|^2  + \frac{4}{r^2} w_2^2 + |w_3'|^2  + \frac{1}{r^2} w_3^2 \right) + \frac{L}{2} \left( |w_4'|^2  + \frac{1}{r^2} w_4^2\right)  \right\} r \, dr \non \\
 &+ \int_0^R \left\{ \left( -\frac{a^2}{2}+\frac{c^2}{4} \left(\sum_{i=0}^4 w_i^2\right) \right) \left( \sum_{i=0}^4 w_i^2 \right) \right.\\
 &\left. - \frac{b^2\sqrt{6}}{36}\left( 2w_0^3-6w_0(w_1^2+w_2^2) +3w_0 (w_3^2+w_4^2) + 3\sqrt{3} w_1(w_3^2-w_4^2) + 6 \sqrt{3} w_2 w_3 w_4 \right) \right\} r \, dr \non
\eea
and using the coercivity conditions $L>0$, $L+\frac43M>0$, we can obtain the existence of a minimizer by  standard methods. 
\eproof
\begin{remark}
It is clear from the system \eqref{sys5} that one can take $w_2=w_3=w_4=0$ and obtain a special $k$-radially symmetric solution of the full Euler-Lagrange equation \eqref{eq:EL} of the form
\be
\label{eq: restricted form}
Q=w_0(r) E_0 + w_1(r)E_1,
\ee
where $(w_0, w_1) \in  S $ minimizes the following energy (see also the corresponding representation \eqref{def:mcR5n})
\bea\label{def:mcR}
\mcE(w_0,w_1)&= \int_0^R \left\{ \frac{L}{2} \left(|w_1'|^2 + |w_0'|^2 + \frac{4}{r^2} w_1^2 \right) + \frac{M}{6}  \left| \sqrt{3} w_1' -w_0' + \frac{2\sqrt{3}}{r} w_1 \right|^2 \right\} r \, dr \nonumber\\
 &+ \int_0^R \left\{ \left( -\frac{a^2}{2}+\frac{c^2}{4} \left(w_0^2+w_1^2 \right) \right) (w_0^2+w_1^2) - \frac{b^2\sqrt{6}}{18}\left( w_0^3-3w_0w_1^2 \right) \right\} r \, dr
\eea
and satisfies the following system of ODEs
 \bea \label{ODEsys}
&\left(L+ \frac{M}{3} \right) \left(w_0''+\frac{w_0'}{r}\right) - \frac{M}{\sqrt{3}} \left( w_1'' + \frac{3 w_1'}{r} \right) = w_0 \left( -a^2 - \frac{b^2}{\sqrt{6}} w_0 + c^2 (w_0^2 + w_1^2) \right) + \frac{ b^2}{\sqrt{6}} w_1^2,  \non \\
& (L+M) \left(w_1''+\frac{w_1'}{r}-\frac{4w_1}{r^2}\right) - \frac{M}{\sqrt{3}} \left( w_0'' - \frac{w_0'}{r} \right)  = w_1 \left( -a^2 + \frac{2 b^2}{\sqrt{6}} w_0 + c^2 (w_0^2 +w_1^2) \right)
\eea
with boundary conditions
\be
 w_0'(0)=0, w_1(0) = 0, \ w_0(R) = -\frac{s_+}{\sqrt{6}}, \  w_1(R)=\frac{s_+}{\sqrt{2}}.
\ee
The solution $Q$ will always have $e_3$ as an eigenvector, and therefore is $2$-radially symmetric according to Definition~\ref{def:k_rad} (see also \cite{DRSZ, INSZ2d1}). 
As shown in Section~\ref{sec3}, numerical calculations indicate that generally $2$-radially symmetric solutions are not always of the restricted form \eqref{eq: restricted form}. 
\end{remark}

We have shown the existence of generally $2$-radially symmetric solutions $Q=\sum_{i=0}^4 w_i(r) E_i$ for any $L>0$ and $M  > -\frac34L$. In general, there is no obvious reason for these solutions not to have all five components present. However,  when $M$ is large enough, we can  show that generally 2-radially symmetric solutions are of a restricted form.

\begin{lemma} \label{lem: large M} Assume $L >0$. For any fixed $M$ denote by $\mcE_M(\{w_i\})$ the energy in \eqref{def:mcR5n}.  For $M \to +\infty$, the following statements hold: 

\begin{enumerate}
\item If $\mcE_M(\{w_i^M\}) \leq C$ with $C$ independent of $M$ then $w_0^M \to  w_0$ weakly in $H^1((0,R); rdr)$,  $w_i^{M}\to  w_i$  weakly in $H^1((0,R); r\,dr)\cap  L^2\left((0,R); \frac{1}{r}\,dr\right)$ for $i=0,\dots,4$ (maybe for a subsequence) and 
\be \label{lincon}
w_2= w_3=0, \ \sqrt{3}(w_1 r^2)'=r^2  w_0' \  \hbox{ a.e. } r\in (0,R).
\ee
\item The energy $\mcE_M$ $\Gamma$-converges in the weak topology of $H^1((0,R); rdr) \times [H^1((0,R); r\,dr)\cap  L^2\left((0,R); \frac{1}{r}\,dr\right)]^4$ to:
$$
\mcE_\infty:=\left\{
\begin{array}{ll} \displaystyle\int_0^R\left\{\frac{L}{2}(|w_1'|^2+|w_0'|^2+|w_4'|^2+\frac{1}{r^2}(4w_1^2+w_4^2))+  h(w_0,w_1,w_4)\right\}r\,dr & \textrm{if }\eqref{lincon} \hbox{ holds}\\
+\infty & \textrm{otherwise,}\end{array}\right.
$$
where
\bea
h(w_0,w_1,w_4) &=  \left( -\frac{a^2}{2}+\frac{c^2}{4} \left( w_0^2 + w_1^2+w_4^2 \right) \right) \left( w_0^2+w_1^2+w_4^2\right) \non \\
 &- \frac{b^2\sqrt{6}}{36}\left( 2w_0^3-6w_0w_1^2 +3w_0 w_4^2  - 3\sqrt{3} w_1w_4^2  \right). \non
\eea
\end{enumerate}
\end{lemma}

\bproof
We notice that the first statement trivially follows from the uniform energy bounds. In order to show $\Gamma$-convergence we need two statements:
 \begin{enumerate}
 \item If $w_0^M \to w_0$ in  $H^1((0,R); rdr)$ and $w_i^M \to w_i$ in $H^1((0,R); r\,dr)\cap  L^2\left((0,R); \frac{1}{r}\,dr\right)$ for $i=1,\dots, 4$  then 
 \bea
 \liminf_{M\to\infty }\mcE_M(\{w_i^M\})\ge \mcE_\infty(\{w_i\}).
 \eea
 \item For any $w_0 \in H^1((0,R); rdr)$ and $w_i \in H^1((0,R); r\,dr)\cap  L^2\left((0,R); \frac{1}{r}\,dr\right)$ for $i =1, \dots, 4$ there exists a sequence $w_0^M \in H^1((0,R); rdr)$ and $w_i^M \in H^1((0,R); r\,dr)\cap  L^2\left((0,R); \frac{1}{r}\,dr\right)$ for $i =1, \dots, 4$  such that
 $$
 \limsup_{M\to\infty}\tilde\mcE_M(w_i^M)=\mcE_\infty(w_i).
 $$
 \end{enumerate}
The first statement is just lower semicontinuity of the energy $\mcE_M$ and for the second we take $w_i^M=w_i,i=0,\dots,4$. The existence of minimizers for the limiting energy $\mcE_\infty$ follows by standard arguments.
\eproof

\section{$k=2$ defects} 
\label{sec3}

In this section, we carry out a numerical study of two-dimensional defects  for the two-elastic-constant Landau-de Gennes energy \eqref{LDG} 
in the case $k = 2$.  From Theorem~\ref{thm:radialsym}, it is only for $k = 2$ that generally $k$-radially symmetric solutions of the Euler-Lagrange equations exist for general nonzero $M$.
 We consider  the cases $b = 0$ and $b=1$ separately in Sections~\ref{subsec: b=0} and \ref{subsec: b=1} respectively. The two cases are distinguished by the set of minimisers of  the bulk potential $f_b(Q)$, which for $b=0$ is
the four-sphere $|Q|^2 =\frac{a^2}{c^2}$, while for $b > 0$ is the two-dimensional space of uniaxial $Q$-tensors with scalar order parameter $s_+$ given by  \eqref{s+}.     The numerical results build upon previous analytic results for the case $M=0$ presented in \cite{DRSZ} for $b=0$ and \cite{INSZ2d1, INSZ2d2} for $b > 0$. Throughout, we fix $L=1$ and $a = c = 1$. 
As explained below, 
 we  find global and local minimisers of the Landau-de Gennes energy as well as some, but not necessarily all, of the generally $2$-radially symmetric saddle points. We find bifurcations in these solutions as functions of $M$ and $R$, where $R$ is the radius of the disk, and study their characteristics  in different parameter regimes.

The numerical procedure consists of finding local minimisers of i) the full
energy  \eqref{LDG},  ii) the reduced energy \eqref{def:mcR5},  which is
defined on  generally $2$-radially symmetric profiles characterised by  radial
functions $w_0(r),\ldots, w_4(r)$, and iii) the restriction of the reduced
energy obtained by taking only $w_0$ and $w_1$ to be nonzero, given by
\eqref{def:mcR}.  From Theorem~\ref{thm:radialsym}, it follows that local
minimisers found in ii) and iii)  are necessarily critical points of the full
energy \eqref{LDG}.  By considering  the second variation of the full energy,
we can determine  whether these critical points are  local minimisers or
saddle points. We also find non-radially-symmetric local minimisers of the
full energy.  Comparing the energies of all the local minimisers obtained in
this way,  we identify the global minimiser as the one with the lowest energy.
As we are not solving the Euler-Lagrange equation directly, we are not able to
find saddle points of the full energy that are not radially symmetric, nor are
we able to find radially symmetric critical points that are not minimisers of
either the reduced energy \eqref{def:mcR5}  or its two-component restriction
\eqref{def:mcR}. We employ a numerical technique of continuation in parameter
$M$ similar to one used in \cite{mg2}, where minimisers of the energy for the
previous value of $M$ are used as initial trials for the perturbed value of
$M$. At each energy minimisation step  for \eqref{LDG} after discretising $Q$-tensors satisfying boundary condition \eqref{BC1} with linear finite elements the resulting
discrete energy is minimised using a trust-region method of nonlinear unconstrained optimisation.
We restrict our numerical analysis to moderate values of the radius, $5 \leq R \leq 50$ and the elastic constant, $0.75< M \leq 100$.

Following this procedure, we find the following types of radially-symmetric and non-radially-symmetric critical points.

\smallskip
\noindent {\it Radially-symmetric:}
\begin{itemize}\item Two-component  solutions $Q_2$ with only $w_0$ and $w_1$ nonzero.  These are the simplest radial solutions, for which $n$, $m$ and $e_3$ are everywhere eigenvectors of the $Q$-tensors.  We identify two sub-types, namely
\begin{itemize}
\item $ Q_2^-$, with $w_0$ and $w_1$ monotonic, $w_0 < 0$ and $w_1 > 0$.  Solutions of this type were studied in  detail in \cite{DRSZ} for $M=b=0$ and in \cite{INSZ2d1, INSZ2d2} for $M=0$, $b \neq 0$.  
\item $Q_2^\pm$, with  $w_0$ changing sign.  In \cite{DRSZ}, it was shown that for $M=b=0$, a global minimizer of the radial energy cannot be of this type.
\end{itemize}  
\item Three-component  solution $Q_3$, with $w_0$, $w_1$ and $w_3$ nonzero.  Neither $n$ nor $e_3$ are everywhere eigenvectors of $Q_3$, although $m$ is.
\item Five-component  solution $Q_5$, with no vanishing components.  The eigenvectors of $Q_5$ do not retain any of the characteristics of the boundary data;  neither  $n$, $m$ nor $e_3$ are eigenvectors of $Q_5$ except at special points (e.g., the origin).
\end{itemize}

\noindent {\it Non-radially-symmetric:} 

A non-radially-symmetric solution belongs to a one-parameter family of solutions generated by rotations about $e_3$.  They are found to have two isolated points where the $Q$-tensor is uniaxial, and about which the $Q$-tensor has the structure of a defect with $k=1$.  The defects are collinear with and equidistant from the origin, and for definiteness may be taken to lie along the $e_1$-axis.  We find  two sub-types, namely 
\begin{itemize}
\item Vertical non-radial solution $\Qnrz$, for which $e_3$ is everywhere an eigenvector.   
The solution inherits  four discrete symmetries from  the boundary data, namely reflection invariance through the $(e_1,e_3)$- and $(e_2,e_3)$-planes, and $\pi$-rotation invariance about the $e_1$- and $e_2$-axes.   
\item Tilted non-radial solution $\Qnr$, for which $e_3$ is not  an eigenvector (except at special points).  The solution inherits  only two discrete symmetries from the boundary data, namely reflection invariance through the $(e_2,e_3)$-plane and $\pi$-rotation invariance about the $e_1$-axis. 
\end{itemize}

\subsection{Case $b=0$}\label{subsec: b=0}

\begin{figure}[ht]
\begin{center}
\hspace*{0.0cm}\includegraphics[width=0.9\textwidth]{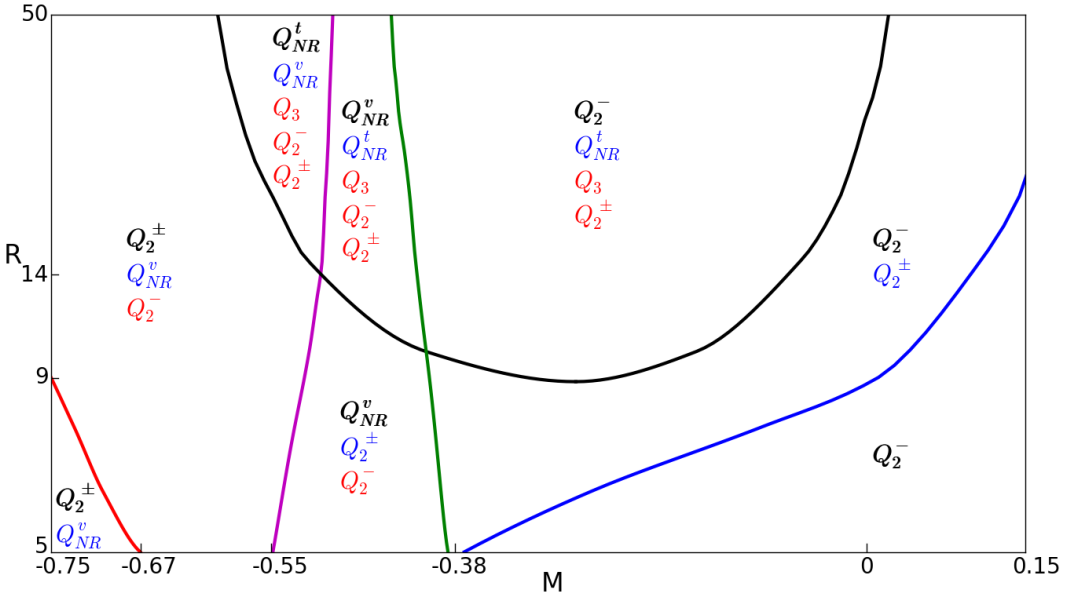}
\caption{\small  Log-log plot of the phase diagram for $b=0$ as a function of elastic constant $M$ and disk radius $R$ scaled as indicated, with $a = c = L =1$. Global minimisers of the Landau-de Gennes energy \eqref{LDG} are indicated in black; local minimisers, in blue; saddle points that are local minimisers of the reduced energy \eqref{def:mcR5} or the 2-component restriction \eqref{def:mcR}, in red.}
\label{fig:b=0_phase_diagram} 
\end{center}
\end{figure}

 Figure~\ref{fig:b=0_phase_diagram} displays a partial phase diagram for the
 Landau-de Gennes energy as a function of the  elastic constant $M$ and
 radius $R$, following the numerical procedure described above.  Global minimisers,
 local minimisers and saddle points are as indicated.  For $M$ sufficiently
 large (to the right of the  blue curve in Figure~\ref{fig:b=0_phase_diagram}), the global minimiser is radially symmetric of type $Q_2^-$, in keeping with  analytic results for  $M=0$ \cite{DRSZ}.  The radially symmetric solution of type $Q_2^\pm$, which  for $M=0$ was shown  to be saddle point for small $L$ \cite{DRSZ}, actually becomes the global minimiser for $M$ sufficiently negative.  For intermediate values of $M$ (between the magenta, black and green curves), the global minimiser is found to be non-radially symmetric of type  $\Qnr$ or $\Qnrz$.
\begin{figure}[ht] 
\includegraphics[width=0.5\textwidth]{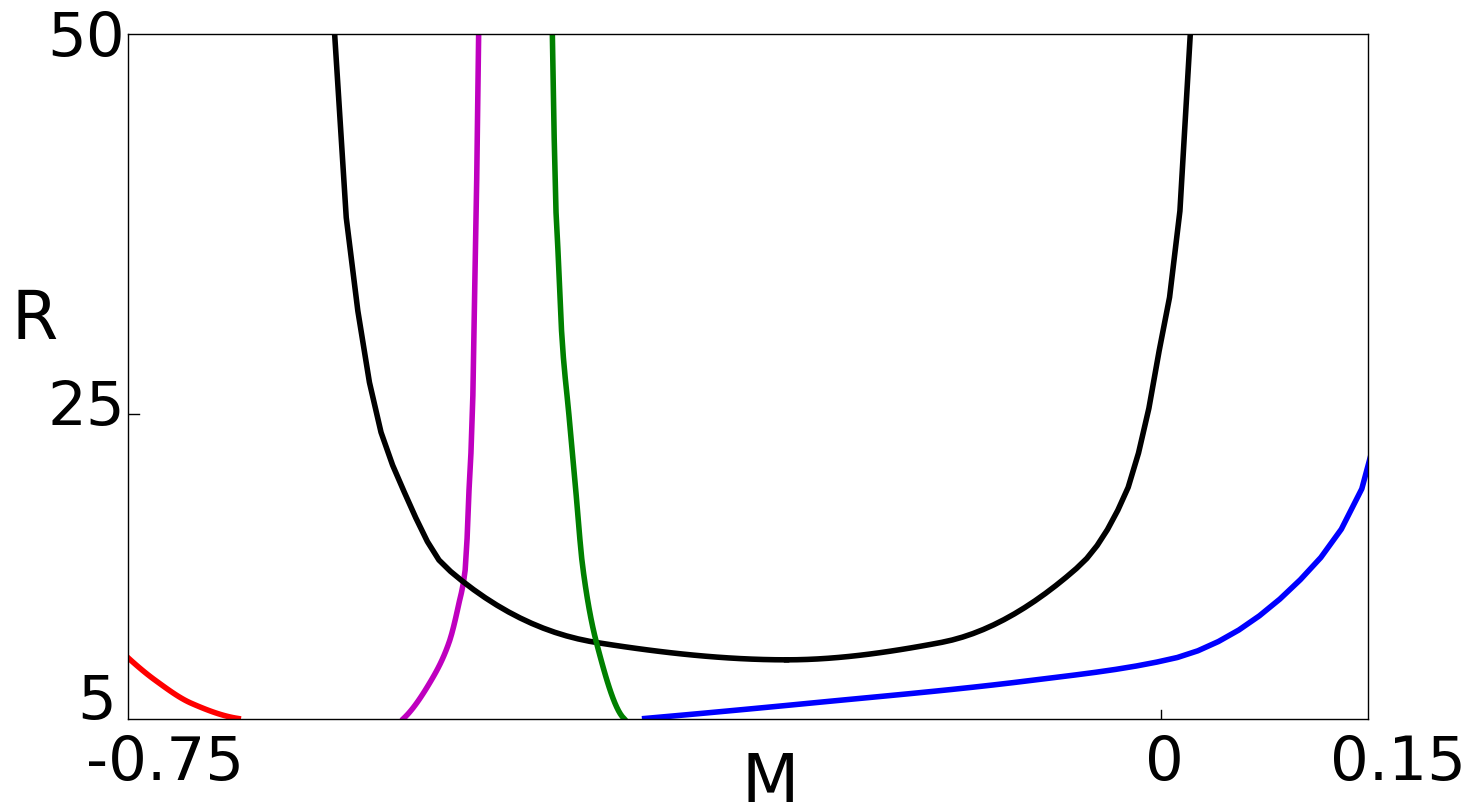}
\includegraphics[width=0.5\textwidth]{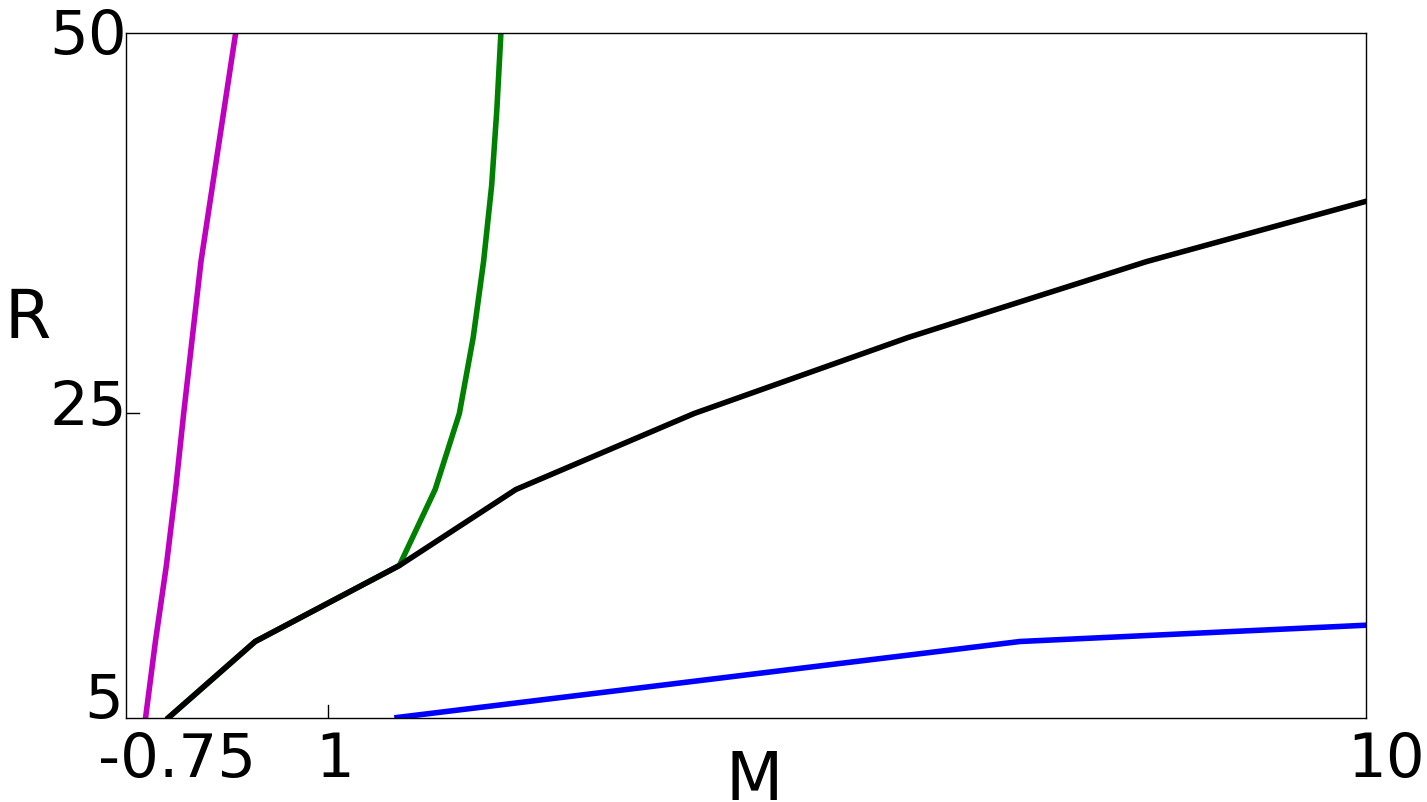}

\caption{\small Unscaled bifurcation diagrams corresponding to Figure~\ref{fig:b=0_phase_diagram} (left) and Figure~\ref{fig:b=1_phase_diagram} (right).}
\label{fig:unscaled_diagrams}
\end{figure}

Even with a partial set of critical points, 
we can discern a rich bifurcation set.   The black curve in
Figure~\ref{fig:b=0_phase_diagram} marks a pitchfork bifurcation, where two three-component radial solutions $Q_3$ (distinguished by the sign of their $w_3$ components) emerge from the two-component solution $Q_2^\pm$, which undergoes a change in stability. 
The fact that $Q_2^-$ is the global minimiser and  $Q_2^\pm$ and $Q_3$  saddle
points  in the region bounded by and above the black curve is consistent with
analytical results  for $M=0$ and small $L$ \cite{DRSZ}. The  black curve also
marks the appearance of the tilted non-radial solution $\Qnr$ as a
 minimiser through a bifurcation with
$Q_2^\pm$.  The blue curve demarcates the domain of existence of $Q_2^\pm$ as a  minimiser of the restricted energy.  The red curve marks the domain of existence of $Q_2^-$ as a  minimiser of the restricted energy.   The magenta curve denotes a first-order phase transition between $Q_2^\pm$ and $\Qnr$ for small $R$ and between $\Qnr$ and $\Qnrz$ for large $R$.

\begin{figure}[H] 
\includegraphics[width=0.35\textwidth]{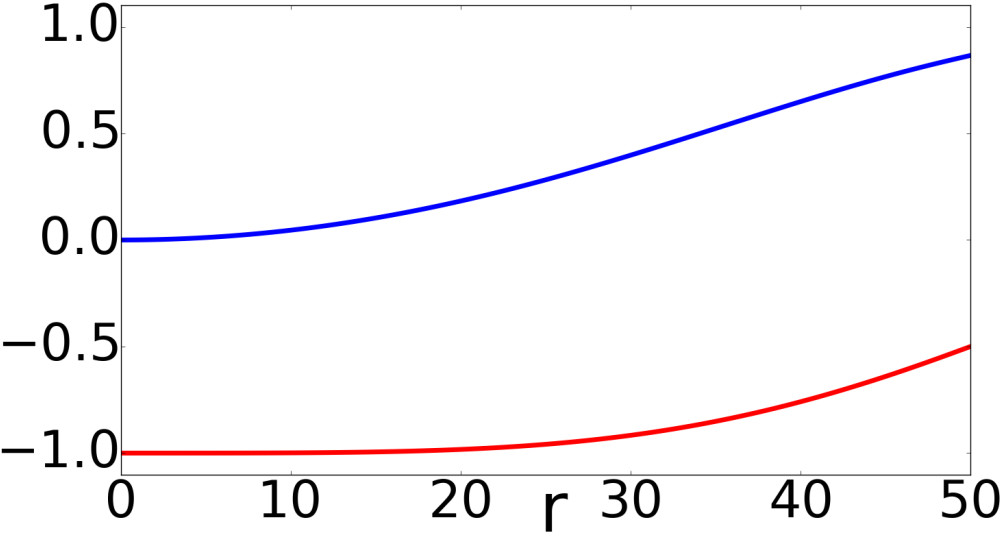}
\hspace{0.01\textwidth}
\includegraphics[width=0.55\textwidth]{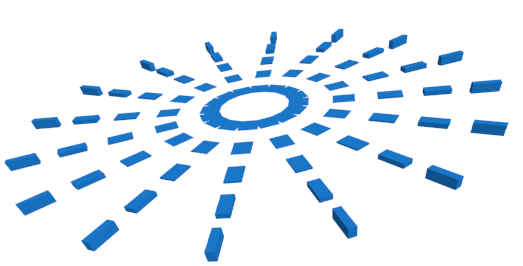}\\[1ex]

\includegraphics[width=0.35\textwidth]{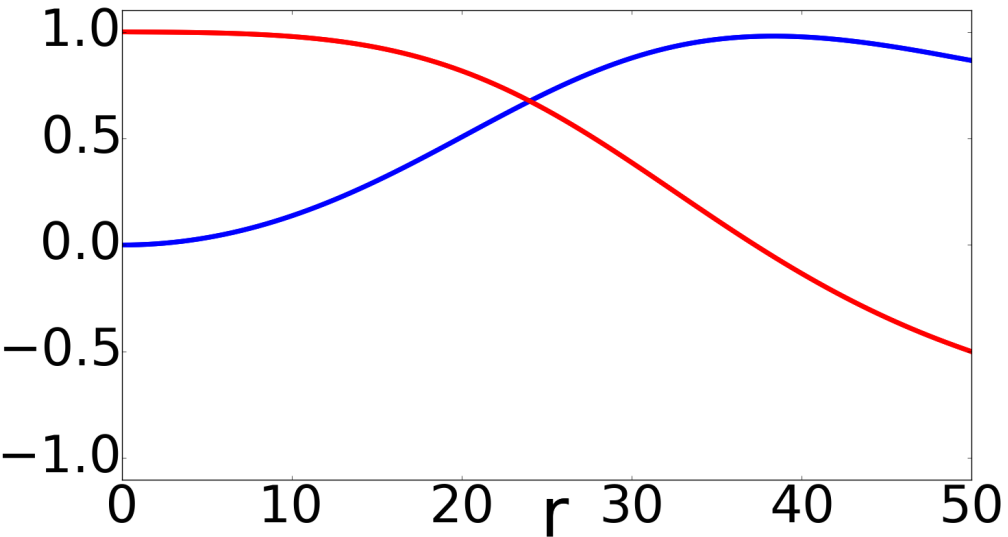}
\hspace{0.01\textwidth}
\includegraphics[width=0.55\textwidth]{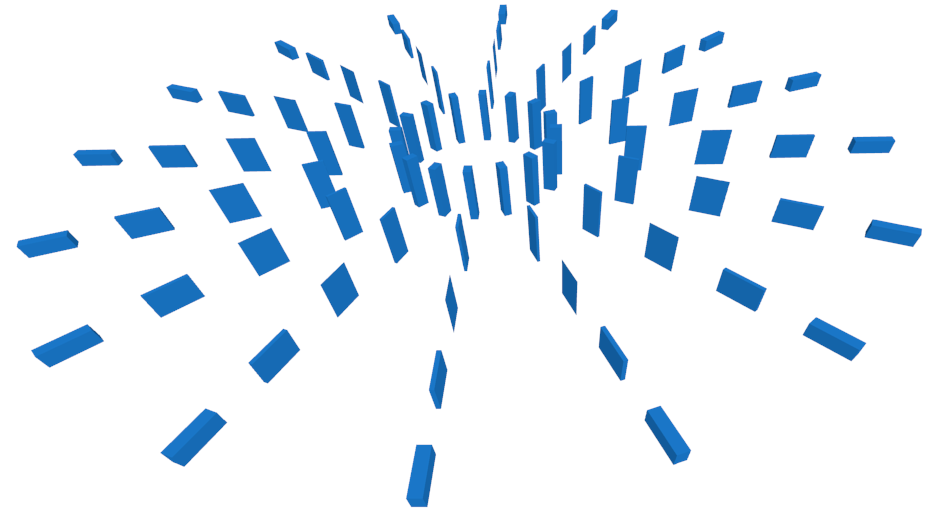}\\[2ex]
\includegraphics[width=0.35\textwidth]{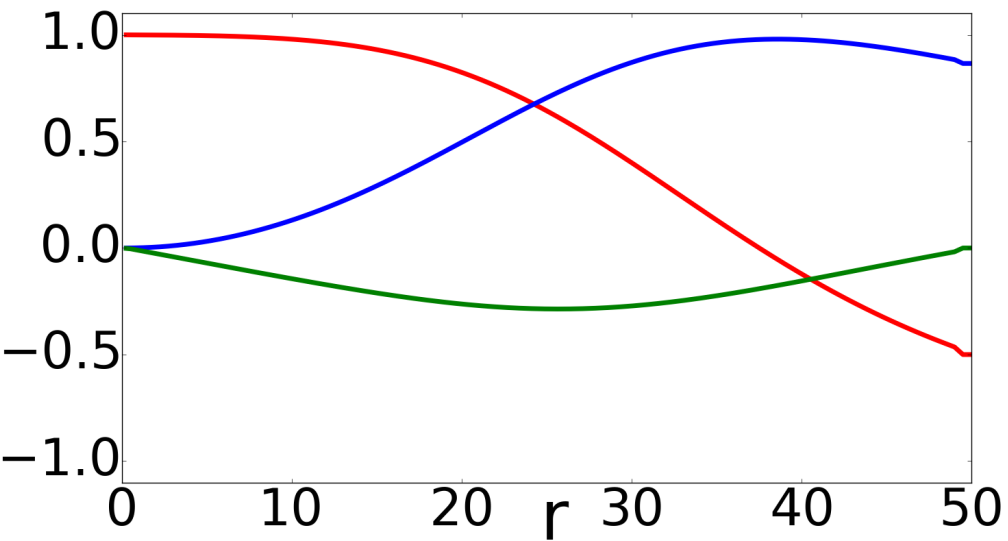}
\hspace{0.01\textwidth}
\includegraphics[width=0.55\textwidth]{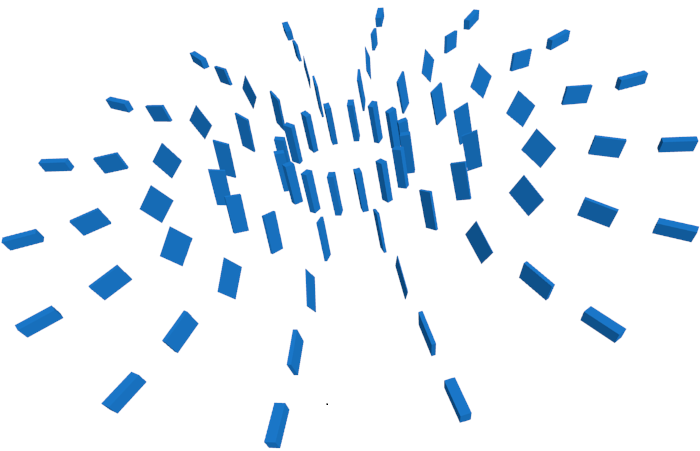}
\caption{\small  Radial components (left) and Q-tensor fields (right) for $Q_2^-$ (top), $Q_2^{\pm}$ (middle), and $Q_3$ (bottom) for $M=0$ and $R = 50$.  As above, $b=0$ and   $a=c=L=1$.
In the radial-component plots on the left, $w_0$ is indicated in red, $w_1$ in blue, and $w_3$ in green.  In the $Q$-tensor plots on the right, $Q$-tensors are represented as parallelepipeds with  axes  parallel to the eigenvectors of  $Q(x)$
and with (nonnegative) lengths given by the eigenvalues of $Q(x)$ augmented by adding $\sqrt\frac23 |Q(x)|$ (see \cite{SonKilHes}).}
\label{2CS_0p1}
\end{figure}

Next, we examine the profiles of the critical points identified in Figure~\ref{fig:b=0_phase_diagram}. The radially-symmetric solutions $Q_2^-$, $Q_2^\pm$ and $Q_3$  are shown in Figure~\ref{2CS_0p1}, including their radial components  and  $Q$-tensor plots. 

We see that the radial components $w_0$ and $w_1$ of $Q_2^-$ are monotonic, 
with $w_0 <  0$ and $w_1 >0$. In contrast, for $Q_2^\pm$,  $w_0$ changes sign  and $w_1$ is not monotonic.  From the $Q$-tensor plots, it is apparent that as the origin is approached, $Q_2^-$ becomes oblate nematic and $Q_2^\pm$ prolate nematic, with the director along  $e_3$ for both.  As noted above, $Q_3$ bifurcates from $Q_2^\pm$  with $w_3$ becoming either positive or negative (the $Q_3$ profile shown in Figure~\ref{2CS_0p1} has $w_3 \le 0$).  The qualitative similarity between the $w_0$ and $w_1$ components of $Q_2^\pm$ and $Q_3$ is apparent.  The $Q$-tensor plot of $Q_3$ shows that the director (ie, the principal eigenvector of $Q_3$) rotates from $e_3$ at the origin to $n$ at the boundary, reminiscent of the escape-to-the-third-dimension (or skyrmion) solution in Oseen-Frank theory \cite{cladiskle}.

The plots shown in Figure~\ref{2CS_0p1} are for $M=0$.  In this case, explicit formulas for the corresponding solutions are provided in the limit $L\rightarrow 0$  by Theorem 4.6 of \cite{DRSZ}, where, in particular, $Q_3$ becomes everywhere uniaxial, and the properties described above may be derived from these limiting formulas.   
For $M\neq 0$, the qualitative behaviour of these radially-symmetric solutions is found to be the same.

\begin{figure}[H] 
\includegraphics[width=0.41\textwidth]{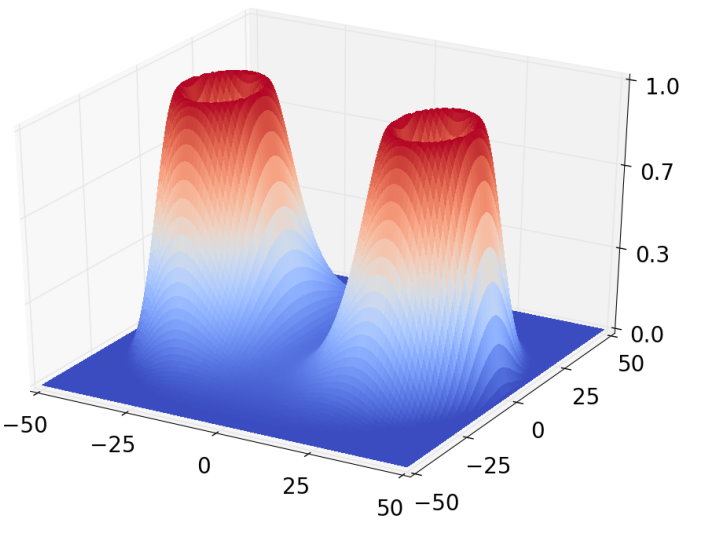}
\includegraphics[width=0.59\textwidth]{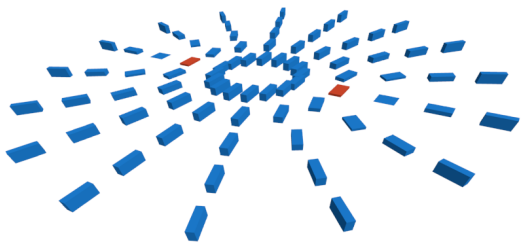}\\[2ex]

\includegraphics[width=0.41\textwidth]{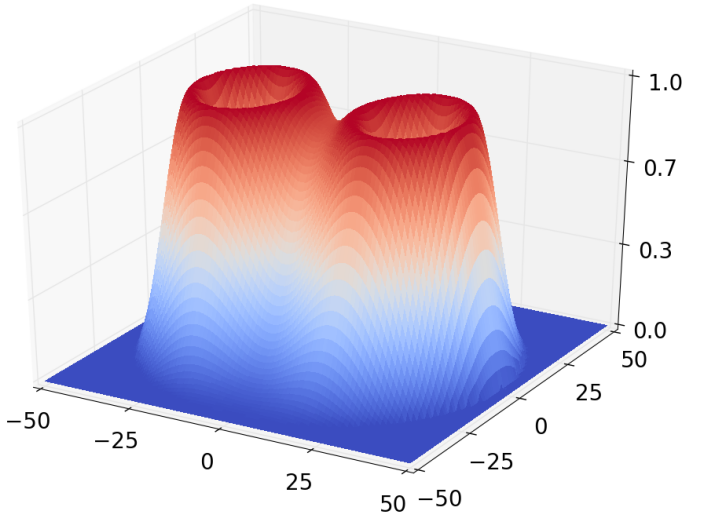}
\includegraphics[width=0.59\textwidth]{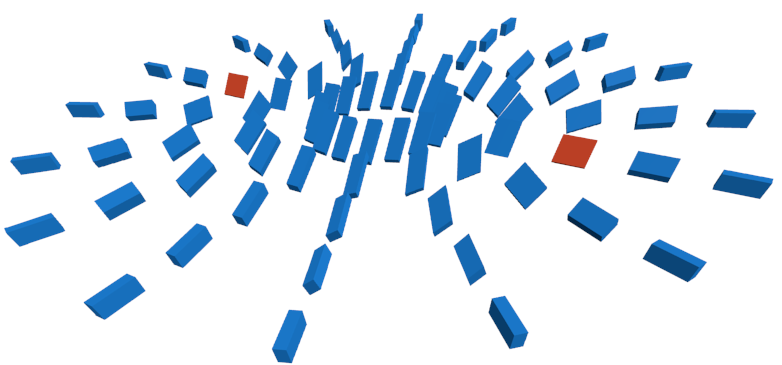}

\caption{\small $\beta$ contours (left) and $Q$-tensor fields (right)  for
  $\Qnrz$ (top) and $\Qnr$ (bottom) for $M=-0.5$ and $R = 50$.  As above,
  $b=0$, and $a=c=L=1$ throughout.  The  location of the defects
  is indicated by the red boxes in the $Q$-tensor plots.}
\label{NSS_0p1}
\end{figure}

The non-radially symmetric solutions $\Qnrz$ and $\Qnr$ are shown in Figure~\ref{NSS_0p1} with  $M = -0.5$ and $R=50$.  We see from 
Figure~\ref{fig:b=0_phase_diagram} that $\Qnrz$ is the global minimiser and $\Qnr$ is a local minimiser in this regime.
In addition to $Q$-tensor fields, we plot the {\it biaxiality parameter}, $\beta$, given by
\be\label{beta}
\beta = 1- 6 \frac{(\tr(Q^3))^2}{(\tr(Q^2))^3}.
\ee
It is straightforward to show that $\beta$ takes values between $0$ and $1$, and that  $\beta =0$ if and only if  $Q$ is uniaxial, while $\beta =1$ if and only if $Q$ is  {\it maximally biaxial}, i.e.~the eigenvalues of $Q$ are equally spaced (equivalently, ${\rm det}\,  Q =0$). The biaxiality plots reveal two half-defects in $\Qnrz$ and $\Qnr$ equally spaced about the origin on the $e_1$-axis, where $\beta = 0$.   
As discussed above, $\Qnrz$ has four discrete symmetries and $\Qnr$ only two.  In particular, $\Qnrz$ has $e_3$ as an eigenvector at the origin, while $\Qnr$ does not, as is apparent in Figure~\ref{NSS_0p1} . 
From the biaxiality plots, it is apparent that neither solution is uniaxial along the $y$-axis, but $\Qnrz$ is much closer to being uniaxial there than is $\Qnr$.

\subsection{Case $b=1$}\label{subsec: b=1}

\begin{figure}[H]
    \centering
           \includegraphics[width=0.9\textwidth]{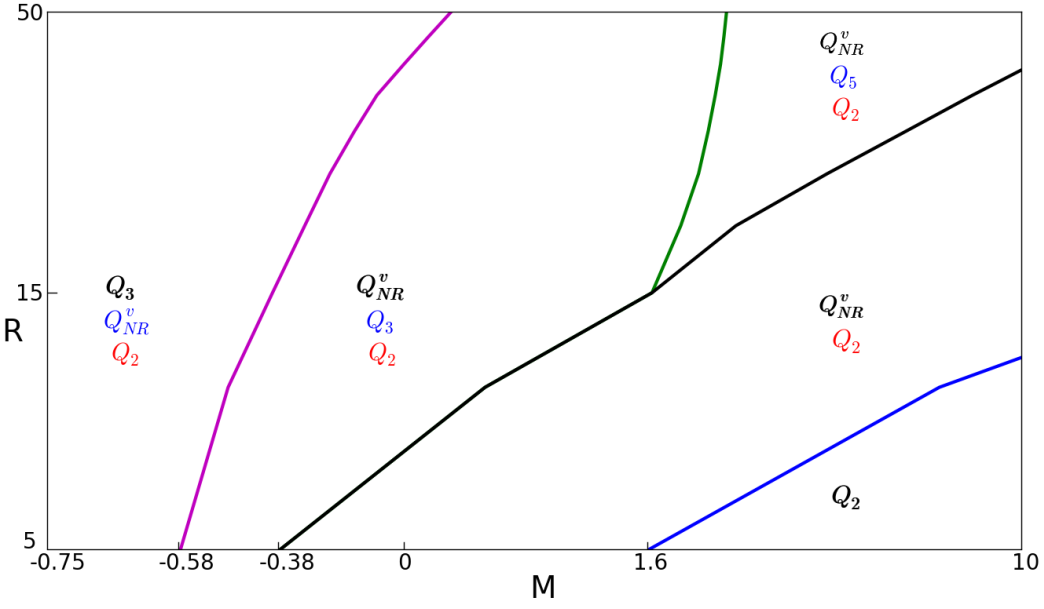}
\caption{\small Log-log plot of existence and stability diagram in case $b=1$.} 
\label{fig:b=1_phase_diagram}
\end{figure}

 Following the numerical procedure described above, we obtain in Figure~\ref{fig:b=1_phase_diagram}  a partial phase diagram for the Landau-de Gennes energy 
with $b = 1$. Global minimisers, local minimisers and saddle points are as
indicated.  We find several differences in the set of solutions in comparison to
the  $b = 0$ case.  Among radially symmetric solutions, we find a single
two-component solution $Q_2$, which is defined everywhere above the coercivity
threshold and which interpolates between the characteristic properties of
$Q_2^-$ 
and $Q_2^\pm$ as $M$ decreases. 
We also find a five-component radial solution $Q_5$,  which was not found for $b=0$.   Among non-radially symmetric solutions, we find only $\Qnrz$. 
 
 In the left-most region of phase plane (large negative $M$), $Q_3$ is the global minimiser.  In the right-most region of the phase plane (sufficiently positive $M$), $Q_2$ is the global minimiser.  In the  region of intermediate $M$, the non-radially symmetric solution $\Qnrz$ is the global minimiser.  For $M=0$, these results are consistent with those of \cite{disclinationsnumerics}.

Even with a partial set of critical points, we find a rich bifurcation set.
The magenta curve describes a first-order
phase transition between $Q_3$ and $\Qnrz$.  The
green curve describes a bifurcation in which $Q_5$ 
emerges from $Q_3$ (note that $Q_3$ is not indicated on the right side of the curve, as it 
ceases to be a local minimiser of the reduced energy
\eqref{def:mcR5} there).  The black curve
describes a bifurcation in which  $Q_5$ and $Q_3$ change stability
for large and small radii respectively. The blue curve describes a bifurcation in which $Q_2$ changes stability and $\Qnrz$ emerges from it through a transcritical bifurcation.  

\begin{figure}[H] 
\includegraphics*[width=0.35\textwidth]{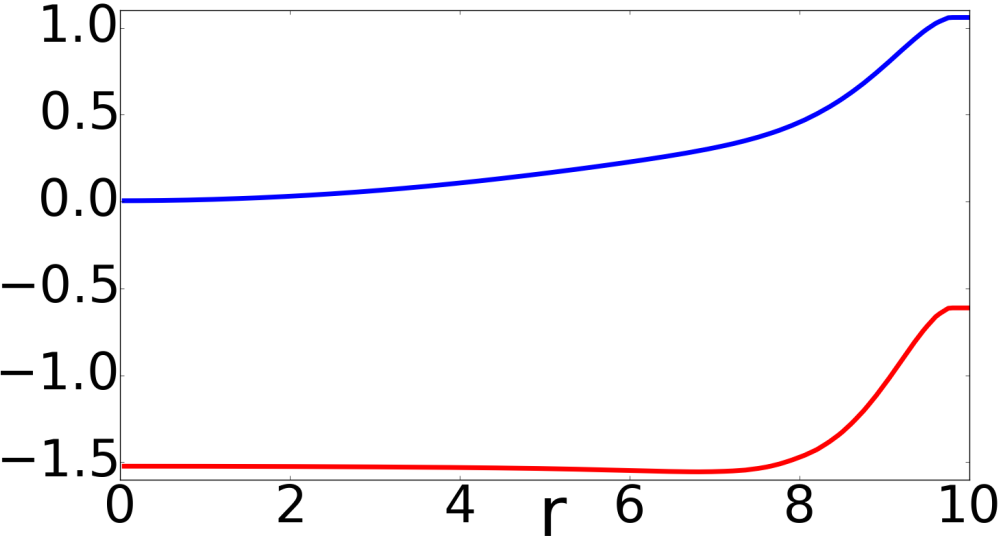}
\hspace{0.01\textwidth}
\includegraphics*[width=0.55\textwidth]{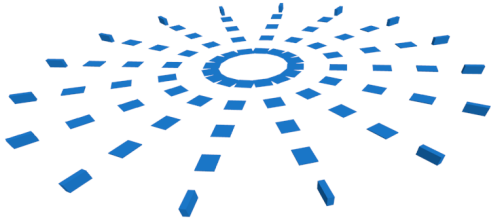}\\[2ex]

\includegraphics*[width=0.35\textwidth]{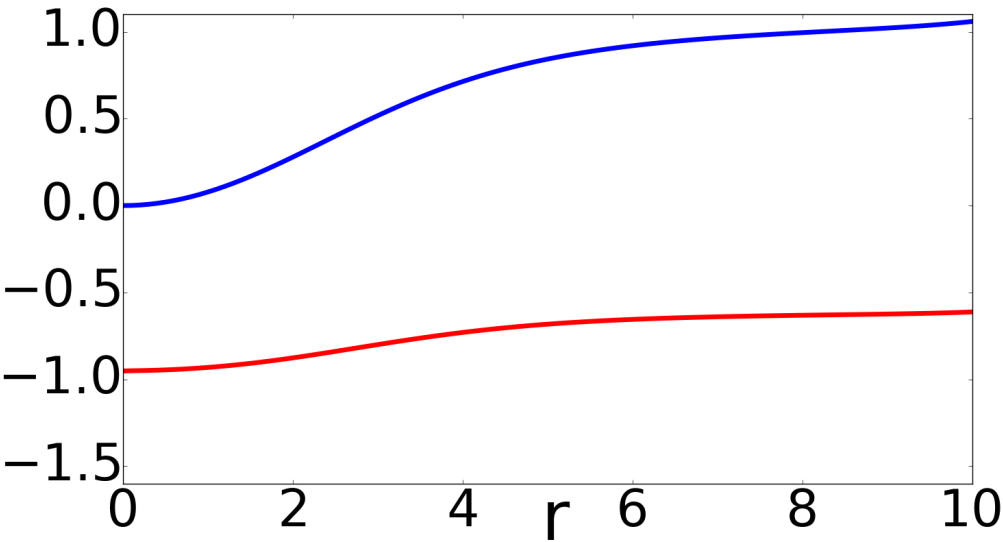}
\includegraphics*[width=0.55\textwidth]{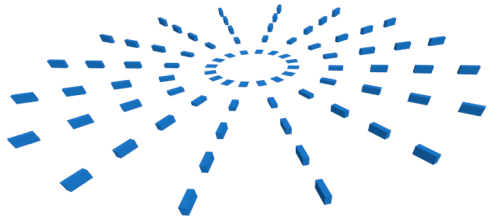}\\[2ex]

\includegraphics*[width=0.35\textwidth]{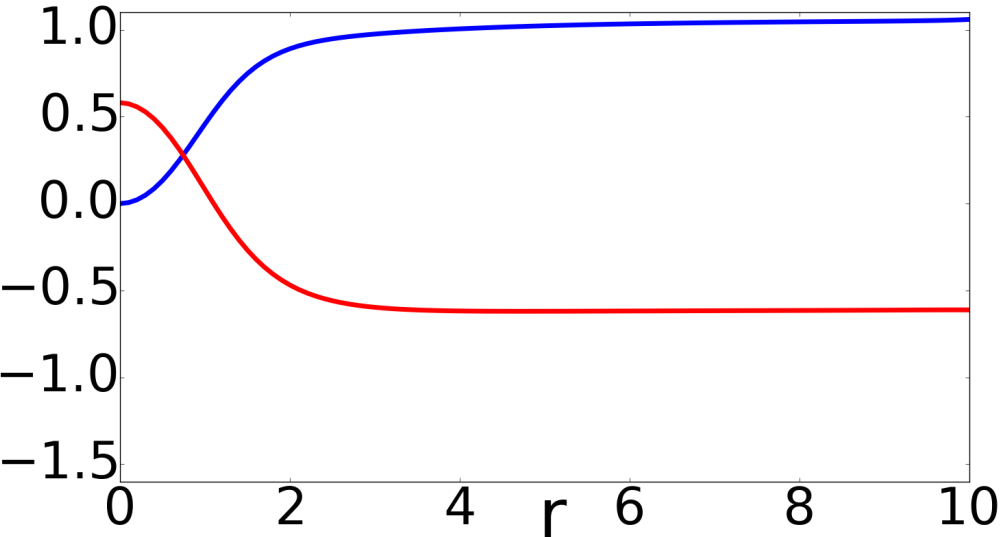}
\includegraphics*[width=0.55\textwidth]{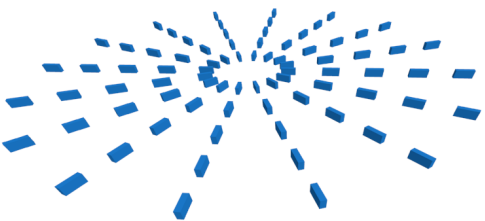}

\caption{\small Radial components (left) and Q-tensor fields (right) of $Q_2$  for $R=10$ and $M=100$ (top row), where $Q_2$ is the global minimiser,  and $M=1$ (middle row) and  $M=-0.55$ (bottom row), where $Q_2$ is a saddle point.  As above, $b=1$ and  $a=c=L=1$ throughout.
In the radial-component plots on the left, $w_0$ is indicated in red and $w_1$ in blue.  In the $Q$-tensor plots on the right, $Q$-tensors are represented as parallelepipeds with  axes  parallel to the eigenvectors of  $Q(x)$ 
and with (nonnegative) lengths given by the eigenvalues of $Q(x)$ augmented by adding $\sqrt\frac23 |Q(x)|$.}
\label{2CS_p}
\end{figure}

Next, we examine the profiles of the critical points identified in Figure~\ref{fig:b=1_phase_diagram}. 
In Figure~\ref{2CS_p} we plot the profiles and $Q$-tensor fields for the two-component solution $Q_2$ for three characteristic  values  of $M$: very large positive $M$,  moderate positive $M$,  and negative $M$. 
In contrast to the  case $b=0$, as $M$ varies and $R$ is kept fixed, the $w_0$ component  of the profile changes continuously from being strictly negative to sign-changing.  As shown in Lemma~\ref{lem: large M}, for large $M$, $w_0$ and $w_1$ are related by the first-order differential equation \eqref{lincon}; for the first profile in Figure~\ref{2CS_p} with $M = 100$, this is verified within numerical accuracy.

\begin{figure}[H] 
\includegraphics*[width=0.39\textwidth]{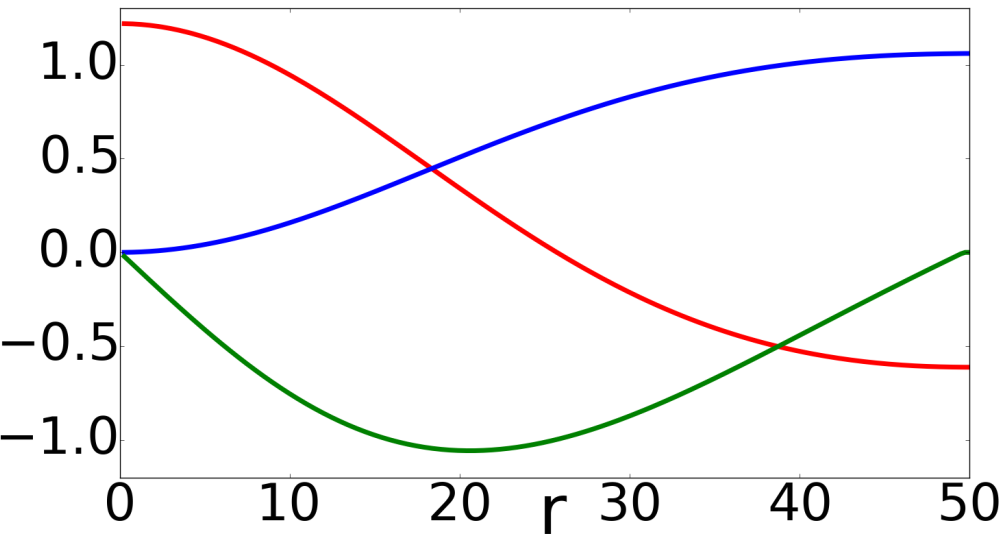}
\includegraphics*[width=0.59\textwidth]{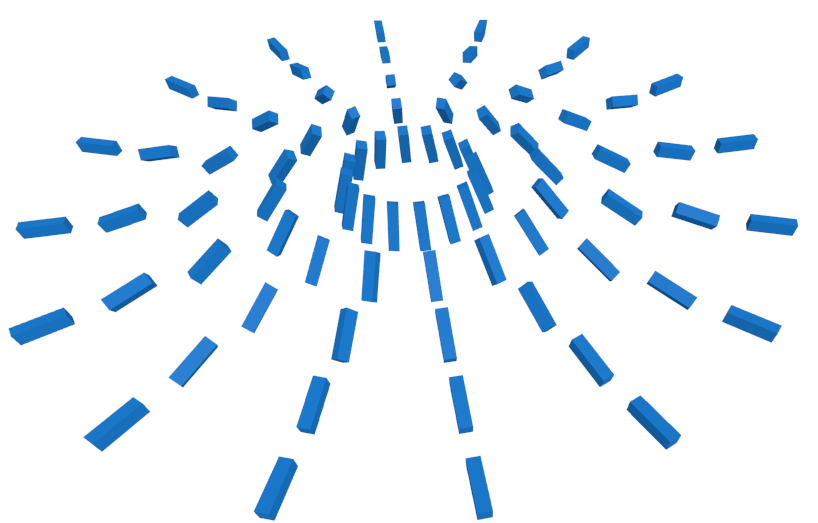}\\[4ex]

\includegraphics*[width=0.39\textwidth]{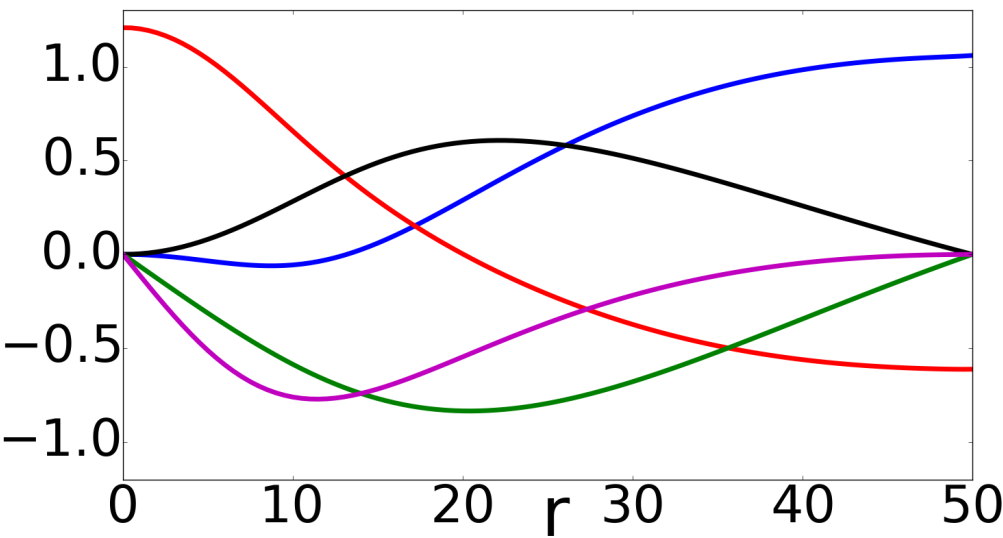}
\includegraphics*[width=0.59\textwidth]{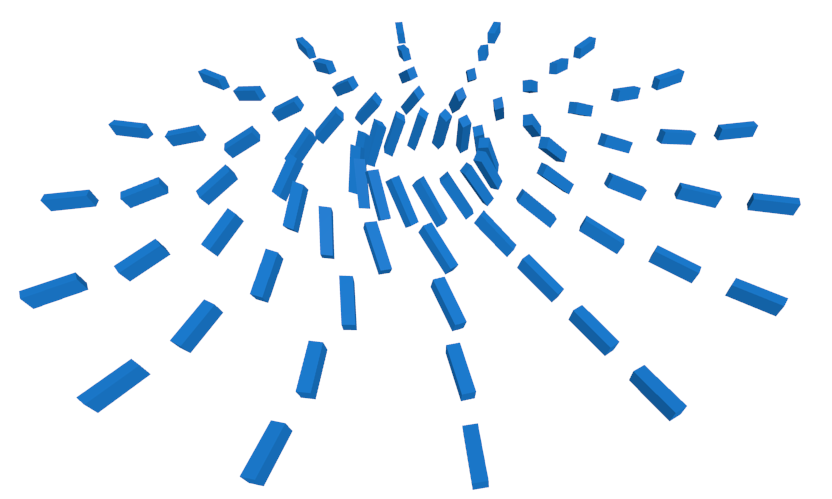}
\caption{\small Radial profiles (left) and Q-tensor fields (right) for $R=50$ for $Q_3$ (top row) with $M = 0$, where $Q_3$ is the global minimiser,  and $Q_5$ (bottom row) with $M = 5$, where $Q_5$ is a local minimiser.  
As above, $b=1$, and  $a=c=L=1$ throughout. 
In the radial-component plots on the left, $w_0$ is indicated in red, $w_1$ in blue, $w_2$ in black, $w_3$ in green and $w_4$ in magenta.}
\label{35CS_p}
\end{figure}

\begin{figure}[H] 
\includegraphics*[width=0.41\textwidth]{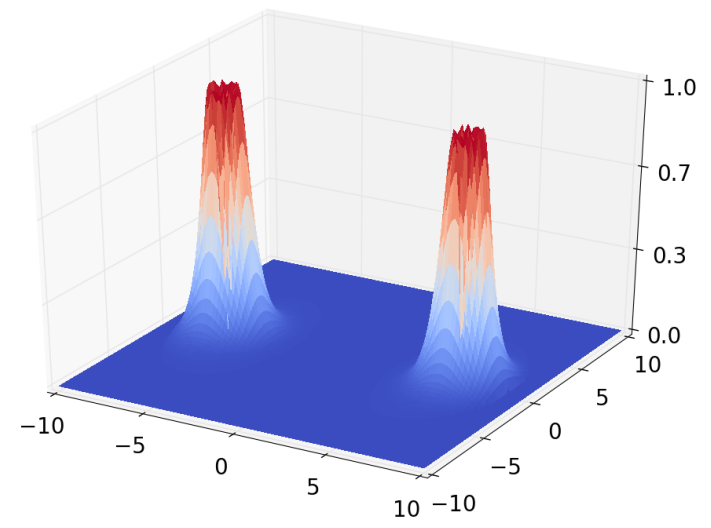}
\includegraphics*[width=0.59\textwidth]{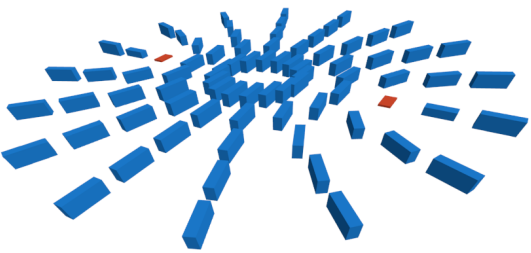}\\[2ex]

\includegraphics*[width=0.41\textwidth]{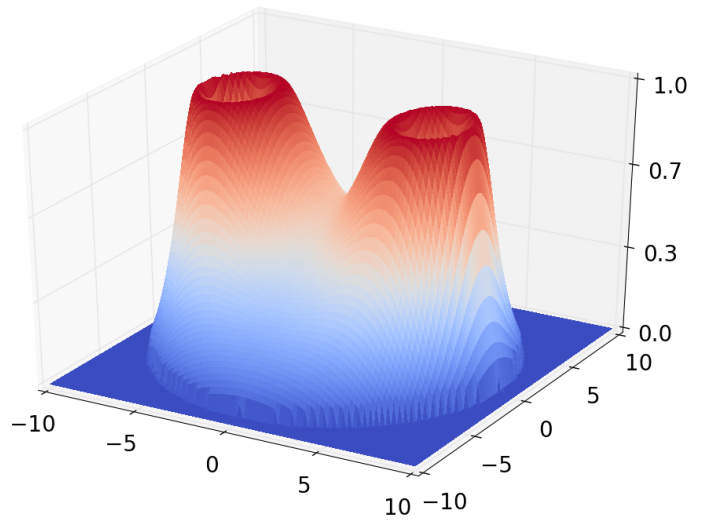}
\includegraphics*[width=0.59\textwidth]{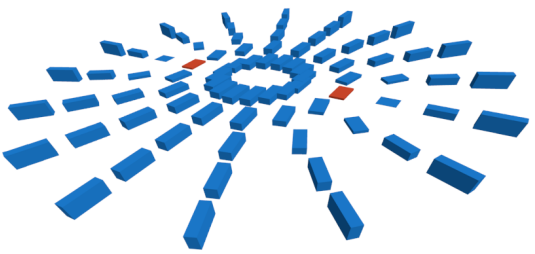}

\caption{\small $\beta$ contours (left) and $Q$-tensor fields (right)  for $\Qnrz$ with $R = 50$ and
 for M=0 (first row) and
  M=5 (second row).  
   As above,
  $b=1$, and $a=c=L=1$.  The  location of the defects
  is indicated by the red boxes in the $Q$-tensor plots.}
\label{NSS_p}
\end{figure}

The radially symmetric solutions $Q_3$ and $Q_5$ are shown in Figure~\ref{35CS_p} for $R=50$ and different values of $M$. 
Recalling that $Q_5$ bifurcates from $Q_3$, we observe that their radial components $w_0$, $w_1$ and $w_3$ are similar.  The $Q$-tensor plots show that $Q_3$ has the structure of the escape-to-the-third-dimension profile \cite{cladiskle}, while $Q_5$ is a twisted version of this profile.
The non-radially symmetric solution $\Qnrz$ is shown in Figure~\ref{NSS_p} for
$R = 10$ and two values of $M$, namely $M=0$ and $M=5$.  From the
$\beta$-contour plots, it is apparent that region of biaxiality is
concentrated around the two defects for $M=0$, while for $M = 5$ it extends to
a neighbourhood of the line joining the defects. 

\section{Symmetry breaking of radially symmetric solutions}
\label{sec4}
\newcommand\mcL{{\mathcal L}}

As shown in Theorem~\ref{thm:radialsym}, for $k\ne 2$ there are no generally radially symmetric critical points of the Landau-de Gennes energy when $M \neq 0$.
We can employ formal perturbation theory  to investigate how  radial symmetry is broken for $M$ nonzero but small.  We restrict our attention to the physically interesting cases $k=\pm 1$.
Here, for $M = 0$, the radial solutions are known to be local minimisers for $b>0$ (see \cite{INSZ2d2}) and global minimisers for $b = 0$ (see \cite{DRSZ}).

Taking $M = \epsilon$, we write the Euler-Lagrange equation in the form
\be L\Delta Q + a^2 Q + b^2 \left( Q^2 - \frac13|Q|^2 I \right) - c^2 Q |Q|^2 = -\epsilon \mcL Q,\ee
where $\mcL Q$ is defined in \eqref{Lop}. 
Let $Y=w_0(r) E_0 + w_1(r) E_1$, with $w_0 < 0$ and $w_1 >0$, denote the radially symmetric solution for the one-constant (i.e., $M = 0$) Landau-de Gennes  energy that is a local minimiser for $b > 0$ \cite{DRSZ} and a global minimiser for $b = 0$  \cite{DRSZ}. Writing $Q = Y + \epsilon W$, substituting above and ignoring terms of order $\epsilon^2$ and higher, we get the following linear inhomogeneous equation for $W$:
 \be \label{eq: eqn for W}  L\Delta W + a^2 W + b^2 \left( YW + WY - \frac23 \tr(YW) I \right) - c^2 W |Y^2| - 2c^2 Y \tr(YW) = - \mcL Y.\ee
The operator on the left-hand side of \eqref{eq: eqn for W} corresponds to the second variation of the one-constant Landau-de Gennes energy evaluated at the radially symmetric solution $Y$.  For $0< b \leq \frac{75}{7} a^2 c^2 $, it has been shown that the second variation at $Y$ is strictly positive definite and it is conjectured that this holds for all $b>0$  \cite{INSZ2d2} .
 (For the special case $b=0$, $Y$ is the  global minimiser, and the second variation is strictly positive for all $k$ \cite{DRSZ}).  Therefore, \eqref{eq: eqn for W} has a unique solution. With some calculation one can show that the solution $W$ is  of the form
\be
\label{a_b_coeff}
W= a_0(r) E_0 + a_1(r) E_1 +  \cos((k-2) \phi) ( b_0(r) E_0 + b_1(r)  E_1)+   \sin((k-2) \phi) b_2(r) E_2,
\ee
where the five radial functions $a_0$, $a_1$ and $b_0$, $b_1$, $b_2$  satisfy a system of linear inhomogeneous ODE's, whose explicit expression has been omitted for brevity. It is clear that $a_0$ and $a_1$ describe perturbations of $w_0$ and $w_1$ respectively, and therefore preserve radial symmetry, while the $b_j$'s comprise the symmetry-breaking component of the perturbed solution, $W_{nr}$, given by
\be \label{eq: Wnr}
W_{nr} =  \cos((k-2) \phi) ( b_0(r) E_0 + b_1(r)  E_1)+   \sin((k-2) \phi) b_2(r) E_2.
\ee

\begin{figure}[H] 
\includegraphics*[width=0.5\textwidth]{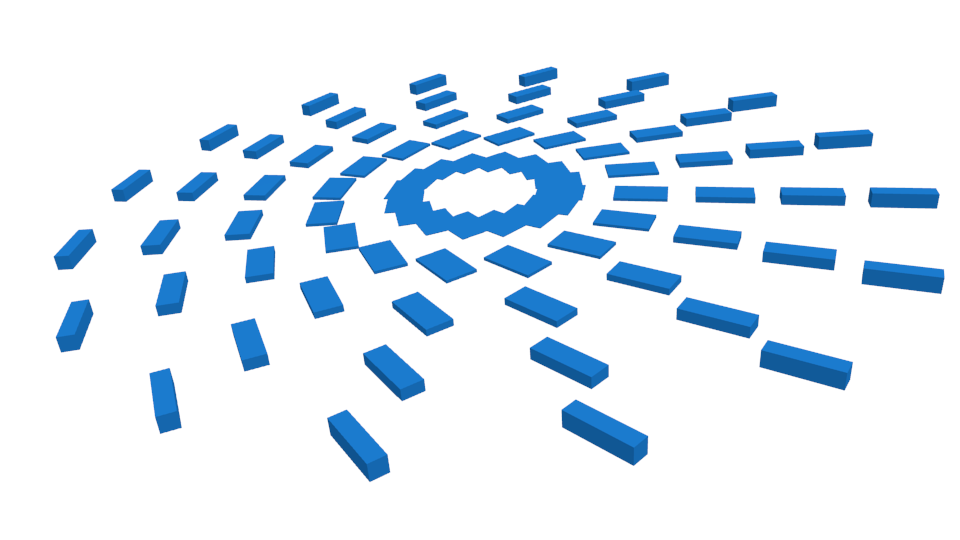}
\includegraphics*[width=0.5\textwidth]{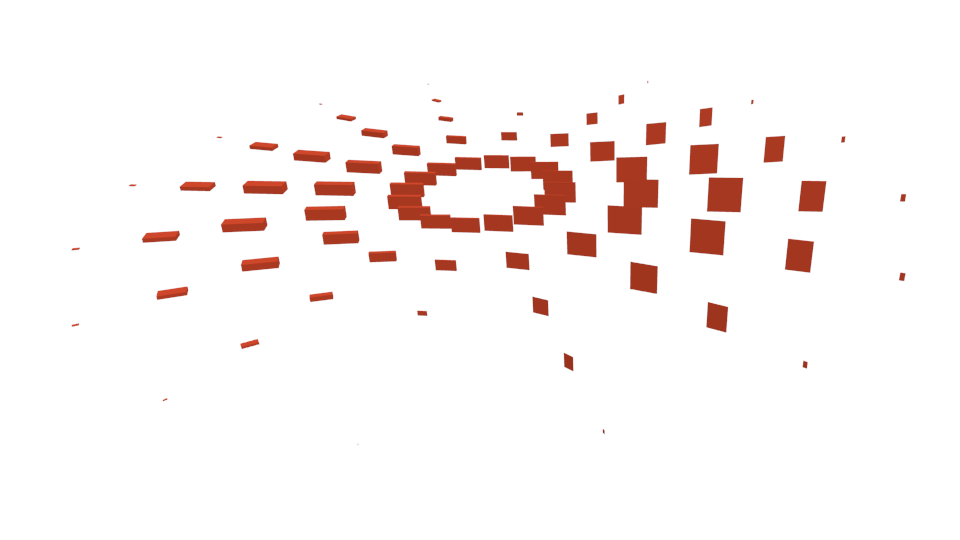}\\[2ex]
\hspace{-0.5cm}
\includegraphics*[width=0.5\textwidth]{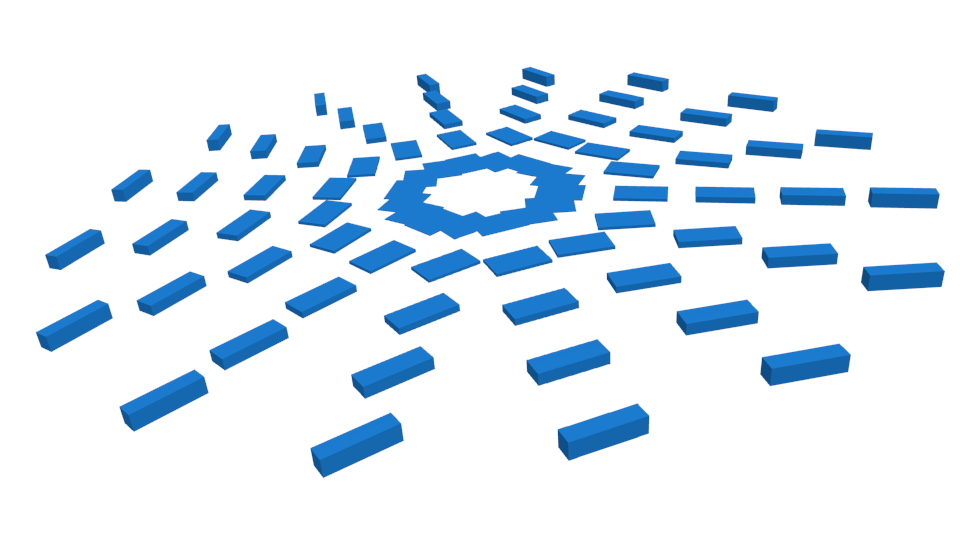}
\includegraphics*[width=0.5\textwidth]{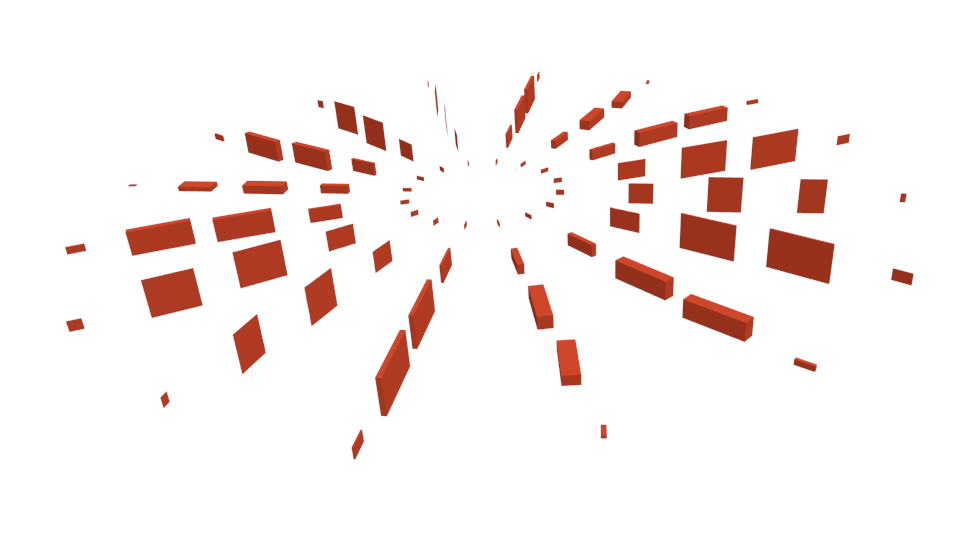}\\[2ex]
\hspace{-0.5cm}
\includegraphics*[width=0.5\textwidth]{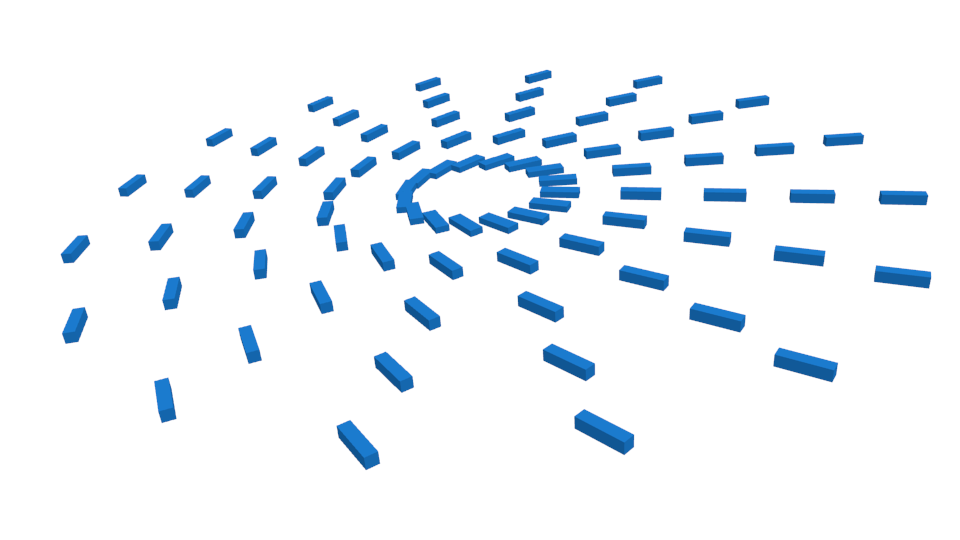}
\includegraphics*[width=0.5\textwidth]{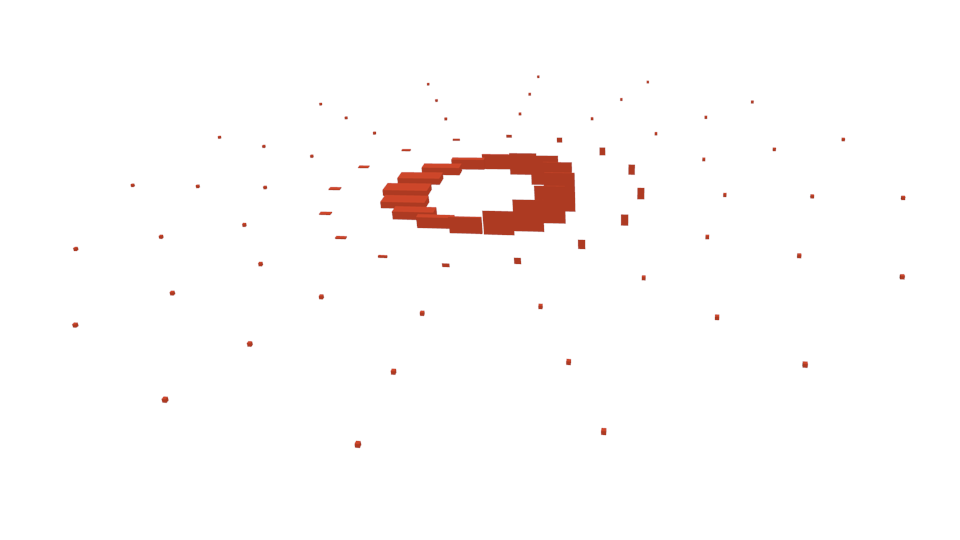}\\[2ex]
\hspace{-0.5cm}
\includegraphics*[width=0.5\textwidth]{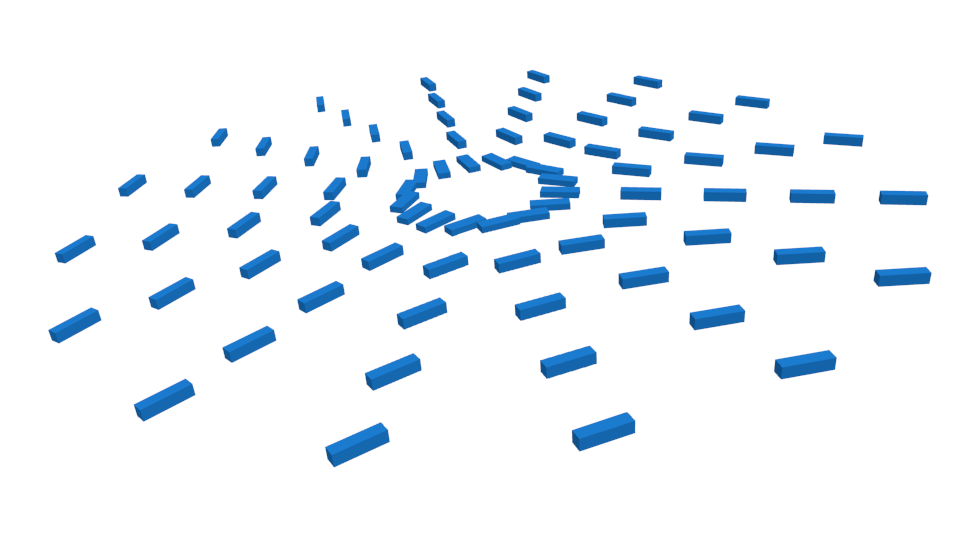}
\includegraphics*[width=0.5\textwidth]{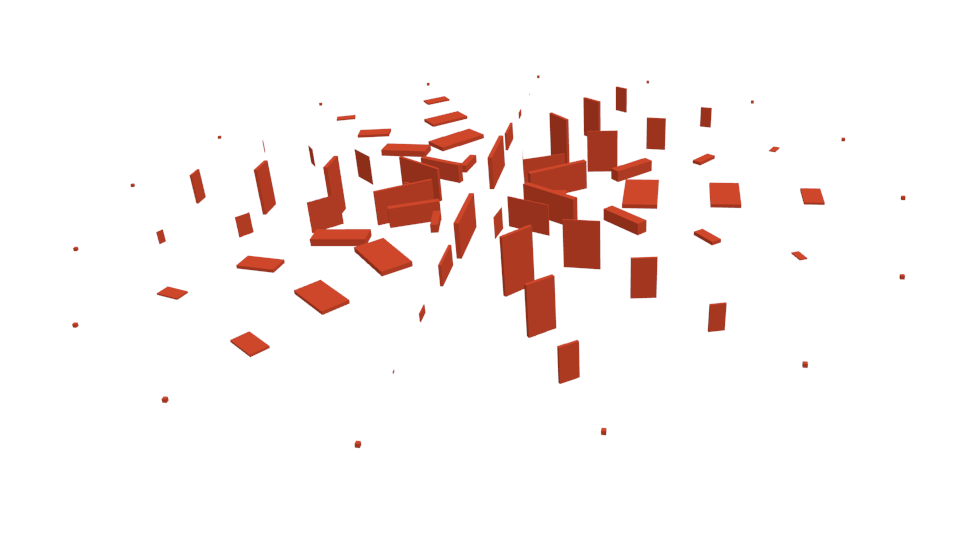}
\caption{\small  $Q$-tensor fields for  the radially symmetric
solutions (left) and the corresponding non-radial perturbations (\ref{eq: Wnr}) (right) for  $b=0,\,k=1,\,R=5$ (first row),
 $b=0,\,k=-1,\,R=5$ (second row),  $b=1,\,k=1,\,R=2$ (third row), and  $b=1,\,k=-1,\,R=2$ (fourth row).
Throughout, $a=20$ and $L=c=1$.  
In the $Q$-tensor fields for the perturbations on the right,  values of $W_{nr}$ very close to zero
are indicated as points.}
\label{Fig:Qtensor}
\end{figure}

It is straightforward  to solve the ODE system numerically to obtain  the $a_j$'s and $b_j$'s.  In  Figure~\ref{Fig:Qtensor},
we give $Q$-tensor plots for the  unperturbed solution $Y$ and symmetry-breaking component $W_{nr}$ for $k=\pm 1$ and $b=0,\,1$.
For $k = -1$, $W_{nr}$ three-fold symmetry, in accord with  \eqref{eq: Wnr}.
The numerical computations confirm that the approximate solutions $Y +
\epsilon W$  are  close to the minimisers $Q_{*\epsilon}$  of the full energy,
and that the norm difference $\Delta := || Q_* - Y - \epsilon W||_{L^2}$
scales as $\epsilon^2$. For example, for the case $k = -1$, $M = \epsilon = 0.1$ and $b =0$, we find that $\Delta/||Y||_{L^2} = 0.024$.   In principle, one can develop a formal perturbation expansion for the exact minimiser in the form $Q_{*\epsilon} = Y + \epsilon W + \epsilon^2 W_2 + \cdots$, and derive a linear inhomogeneous system of ODEs for $W_{n+1}$ in terms of  $Y$ and $W,\ldots, W_n$.   

The computations also reveal that the symmetry between $k$ and $-k$ is broken for $M\neq 0$.  
We note that for $M = 0$,  the  functions $w_0(r)$ and $w_1(r)$ for $k = 1$ and $k = -1$ are the same; the Euler-Lagrange equations depend only on $k^2$ (the $Q$-tensor configurations for $k=1$ and $k=-1$ are not the same, of course, as $E_1$ depends on the sign of $k$).  For $M \neq 0$,
the $k\rightarrow -k$ symmetry is broken, as is readily seen in Figure~\ref{Fig:Qtensor}.  This symmetry breaking in $k$ appears only in  $W_{nr}$, as the $a_0$ and $a_1$ components depend only on $k^2$.

\newpage
\vspace*{-16ex}


\section{Conclusions and open problems }
\label{sec5}

The Landau-de Gennes theory with several elastic constants describes an extended range of rich behaviours in nematic liquid crystals as compared to the single-constant theory. Here we have considered several basic problems related to the existence and behaviour of generally radially symmetric critical points and their symmetry breaking. We have shown that generally radially symmetric  critical points exist only in the case $k=2$, and we have investigated numerically the profiles and stability of these solutions. In particular, we have identified three types of radial profiles: 

\smallskip

(a) {\it two-component} critical point $Q_2=w_0(r) E_0 + w_1(r) E_1$; 

\smallskip 

(b) {\it three-component} critical point $Q_3=w_0(r) E_0 + w_1(r) E_1 + w_3(r) E_3$; 

\smallskip 
(c) {\it five-component} critical point $Q_5=\sum_{i=0}^4 w_i(r) E_i$

\bigskip\par\noindent and we have numerically investigated their minimality depending on the parameters $M$, $R$, and $b$ alongside non-symmetric solutions. Using the minimality property of two-component radial solutions in the case $M=0$, $k=\pm1$ we have studied the symmetry breaking of these minimizers for small $M \neq 0$, identified directions of symmetry-breaking as solutions of a system of ODEs, and verified the principal results of this analysis through numerical simulations.

\smallskip
The model problem we have studied provides an excellent starting point for further analytical investigation of defects in the Landau-de Gennes theory. It is perhaps the simplest possible setting for developing the analytical tools necessary for attacking more general problems, particularly those in three-dimensional domains. Our investigations point  to several natural  open problems whose  solution should lead to significant advances in the area.

\vskip 0.2cm

%
%
%
\noindent {\bf 1. Profiles of radial solutions for $k=2$}

\smallskip
\noindent  There are several important problems related to the qualitative features of the profiles of the radial solutions obtained in the Theorem~\ref{thm:radialsym}. In the case $M=0$ the components $w_0$ and $w_1$ of $2$-radial solution $Q_2^-$ have constant sign and are monotone (see \cite{DRSZ, INSZ2d1, INSZ2d2} for analysis in 2D and \cite{ODE_INSZ} for analysis of a related 3D problem). This fact was particularly useful in the stability analysis when employing the Hardy decomposition trick (see \cite{DRSZ, INSZ2d1, INSZ2d2}). 

\smallskip

\noindent {\bf (a)} Our numerical simulations indicate that for $M>0$ and $b>0$, the components $w_0$ and $w_1$ of the two-component solution $Q=w_0(r) E_0 + w_1(r) E_1$ have definite sign. However, for $M<0$ we observe that $w_0$  changes  sign. It is an interesting problem to explain this behaviour analytically. 

\smallskip 

\noindent {\bf (b)} We observe numerically that for $M<0$ and small $b^2$, the radial components $w_0$ and $w_1$ of the two-component solution are no longer monotonic. It is an interesting problem to study the  monotonicity of these radial profiles. Proving non-monotonicity is a challenging task, since there are no standard techniques available. This problem is also an opportunity to investigate the Landau-de Gennes system in a regime where it is no longer diagonal in the derivative terms. 

\smallskip
\noindent {\bf (c)} 
%
A challenging problem is to establish the existence and stability of the three-component and five-component radial solutions $Q_3$ and $Q_5$.
We note that in some regimes (see Figures~\ref{2CS_0p1} and \ref{35CS_p}), the  components of these solutions exhibit certain characteristics (eg,  sign definiteness, monotonicity, convexity). Investigating these qualitative features may help to understand the stability properties of these solutions.

\bigskip

\noindent { \bf 2. Multiplicity of the local minimizers and critical points}

\smallskip
\noindent  The problem of determining exactly, or even  just lower and upper bounds for, the number  of solutions for the nonlinear system of PDEs (or ODEs) is a very challenging and interesting problem.   There exist several methods of nonlinear analysis, usually related to degree theory, that allow to show the existence of several critical points/local minimizers. However identifying {\it all} critical points (or even local minimizers) requires additional information, specific to the problem. The problems identified in this paper combine low dimensionality with the complexity of a coupled system of PDEs (or ODEs), and present an ideal playground for studying the question of multiplicity of solutions. 

\smallskip

\noindent {\bf (a)} Numerical minimization of the energy \eqref{def:mcR} indicates existence of at least three solutions in some range of parameters $M$ and $R$. Finding the lower and upper bounds on the number of critical points of this energy is the simplest possible problem related to the multiplicity of critical points of the full Landau-de Gennes energy. But even this problem presents a significant challenge due to nonlinear coupling in the Euler-Lagrange equations. We suggest to separately study cases $b=0$ and $b \neq 0$, as in the former case the equations are partially decoupled due to special form of nonlinearity.

\smallskip

\noindent {\bf (b)} Finding  lower and upper bounds on the number of critical points of the full radial energy \eqref{def:mcR5n} is a challenging mathematical problem. Numerical simulations indicate strong dependence of the number of critical points (and local minimizers) on $M$, $R$, and $b$. 

\smallskip

\noindent {\bf (c)} Identifying the exact number of critical point of the full Euler-Lagrange equations \eqref{eq:EL}, \eqref{BC1} seems out of reach by  current analytical methods. However, identifying  lower and upper bounds on the number of symmetry-breaking solutions might be a first step towards the ultimate goal of classifying critical points of the Landau-de Gennes energy. For $M=0$, numerical simulations for $b>0$ indicate that there is at least one non-radial solution with ${\mathbb Z}_k$-symmetry \cite{disclinationsnumerics}.

\bigskip
\noindent {\bf 3. Bifurcations between various solutions}

\smallskip

Reconstruction of the bifurcation diagram is a question related to the multiplicity of solutions but involves more refined information about the system. Bifurcation problems are interesting and tractable because there are a number of powerful tools for studying them, see for instance \cite{kielhofer}. A crucial step in studying various bifurcations in our model is the linearization of the Euler-Lagrange  system near the known solution. The linear operator obtained in this way is essentially  the quadratic form given by the second variation of the Landau-de Gennes energy. Thus, there exists a natural connection between bifurcation problems and the stability studies previously undertaken for these solutions \cite{DRSZ, INSZ2d1, INSZ2d2}. Below are some bifurcation problems inspired by the numerical explorations reported here.

\bigskip 
\noindent {\bf (a)} 
There are three natural parameters, namely the radius $R$, the bulk parameter $b$ and the elastic constant $M$, and one can attempt to study bifurcations with respect to any one of these, holding the others fixed. 
One possible way to understand  symmetry-breaking solutions is to use the fact near the bifurcation, the symmetry-breaking solution will be close to a radially symmetric solution, and the difference between them is related to the kernel of the linearisation at the radially symmetric solution (see for instance the standard bifurcation theorems of Crandall-Rabinowitz type in \cite{kielhofer}).

\smallskip

\noindent {\bf (b)}
In the case $k=2,$ there exist several types of radially symmetric solutions, depending on the value of $M$.  A  study of bifurcations of these radially symmetric solutions with  $M$ as bifurcation parameter appears to be a simpler problem, but should provide important information about various types of critical points of the Landau-de Gennes energy.

\bigskip

\noindent {\bf 4. Global minimality } 

\smallskip

\noindent Identifying a global minimizer is one of the most challenging problems in the calculus of variations. Most of the current methods allow for finding local minimizers, while proofs of global minimality are scarce and problem-dependent. Below we list several specific problems where we believe it may be possible to identify a global minimizer.

\smallskip

\noindent {\bf (a)} It was shown in \cite{DRSZ} that for $b=0$ and $M=0$, the 2-radially symmetric solution is the global minimizer of the Landau-de Gennes energy. It is an interesting question to determine to what extent this statement remains true in the case $M\not=0$. We observe numerically that the two-component radial solution is the global minimizer for a range of $M$ and $R$ (see Figure~\ref{fig:b=0_phase_diagram}). However, the qualitative features of the components (in particular being sign changing) depend on $M$ and $R$. Therefore it is not straightforward to adapt the proof in \cite{DRSZ} for these regimes. A new proof of global minimality should probably depend on $M$ and $R$. 

\smallskip
\noindent {\bf (b)} In the case $M=0$ and $k =\pm1$,  the results in \cite{INSZ2d2} show that the two-component radially symmetric solution $Q=w_0(r) E_0 + w_1(r) E_1$  is a local minimizer. Our numerical simulations (see also \cite{disclinationsnumerics}) suggest that this solutions is a global minimizer of Landau-de Gennes energy. It is an interesting and challenging problem to prove this fact analytically. Related but even more difficult is the question of uniqueness: are there  any other local minimizers (or critical points) of the Landau-de Gennes energy \eqref{LDG} with boundary conditions \eqref{BC1} for $M=0$, $k=\pm1$.

%

\section*{Acknowledgment.} 
GK, JMR, VS would like to acknowledge support from Leverhulme grant RPG-2014-226. JMR and VS acknowledge support from EPSRC grant EP/K02390X/1. AZ was partially supported by a grant of the Romanian National Authority for Scientific Research and Innovation, CNCS-UEFISCDI, project number PN-II-RU-TE-2014-4-0657. AZ gratefully acknowledges the hospitality of the Mathematics Department at the University of Bristol, through Leverhulme grant RPG-2014-226. VS and AZ gratefully acknowledge stimulating discussions with Prof. Radu Ignat and
Prof. Luc Nguyen.

\section{Appendix}
\label{appendix}

\subsection{Null Lagrangian}
In this section we show that elastic terms with $L_2$ and $L_3$ are equivalent up to a boundary term. We prove the result in $\RR^3$ and state the analogous result in $\RR^2$ .

\bigskip
\noindent {\bf Proof of Lemma~\ref{lemma:nulllagrangian}.} 
We assume $\Omega \subset \RR^3$ with $C^1$ regular boundary  $\partial \Omega$  and denote by $n$ the exterior unit normal to $\partial\Omega$. 
%
%
We define 
\be
\mcI_2\defeq \int_\Omega Q_{ik,j} Q_{ij,k},\quad \mcI_3\defeq \int_\Omega Q_{ij,j}Q_{ik,k} , \quad 1 \leq i,j,k \leq 3.
\ee
 Straightforward calculation gives
\bea
\mcI_2-\mcI_3&=\int_\Omega (Q_{ij}Q_{ik,j})_{,k}-(Q_{ij}Q_{ik,k})_{,j}\,dx=\int_{\partial\Omega} Q_{ij}Q_{ik,j}n_k-Q_{ij}Q_{ik,k}n_j\,d\sigma \non\\
&=\int_{\partial\Omega}Q_{ij}Q_{il,r}(\delta_{jr} - n_jn_r) n_l-Q_{ij}Q_{ik,l}\left(\delta_{kl}-n_kn_l\right)n_j\,d\sigma \non\\
&=\int_{\partial\Omega} W_{rj} T_{rj} \,d\sigma, \label{claim}
\eea
where
\be W_{rj} = n_l (Q_{,r} Q - Q Q_{,r})_{lj} \ee
and $T := I - n\times n$ is the orthogonal projection onto the tangent space of $\partial \Omega$.
It is evident that $W_{rj} T_{rj}$ involves only $Q$ and its tangential derivatives on $\partial \Omega$, so that with
%
prescribed Dirichlet boundary data  $Q(x) = Q_b (x)$ for all $x \in \partial \Omega$,  the value of $\mcI_2-\mcI_3$ is determined by $Q_b$.  

In order to obtain a relation between $\mcI_2$ and $\mcI_3$ in $\RR^2$ we just assume that $Q$ is independent of variable $x_3$. Using \eqref{claim} it is straightforward to obtain the analogous result for  $\Omega \subset \RR^2$ ,
\be\label{claim1}
\mcI_2-\mcI_3=\int_{\partial\Omega} t_j \left( \frac{\partial Q_{il}}{\partial t}Q_{ij} - \frac{\partial Q_{ij}}{\partial t}Q_{il} \right) n_l \,d\sigma, 
\ee 
where $n =(\hat n, 0)$, $t = (\hat t, 0)$, and $\hat n = (n_1,n_2)$, $\hat t =(t_1, t_2)$  are normal and tangent to $\partial \Omega$, respectively. 
\eproof

\subsection{Coercivity conditions}


In this section we want to derive the coercivity conditions for the energy \eqref{LDG} in $\RR^3$ and $\RR^2$. The elastic part of the energy is
\begin{align}\label{LDG-el}
\mcF_{el} [Q; \Omega]= \int_{\Omega} \Big[\frac{L_1}{2}|\nabla{Q}|^2 +\frac{L_2}{2}\partial_j Q_{ik}\partial_k Q_{ij}+\frac{L_3}{2}\partial_j Q_{ij}\partial_k Q_{ik} \Big]\,dx.
\end{align} 
We prove Lemma~\ref{lemma:coercivity} providing the necessary and sufficient relations between the elastic constants in order to have the coercivity conditions.

\bigskip
\noindent{\bf Proof of Lemma~\ref{lemma:coercivity}.} In the following we argue for smooth $Q$-tensors, the reduction to the $H^1$ case being standard. As usual, we assume throughout summation over repeated indices.

\bigskip
\noindent {\bf 1.} We denote by $P_{ijk}, i,j,k=1,2,3$ the  components of a third-order tensor 
$$
P \in \mcP=\{ P \in \RR^{3 \times 3\times 3}\, :\, P_{ijk}=P_{jik}, \ P_{iik} =0 \hbox{ for all } 1 \leq i, j, k \leq 3\} 
$$ 
and  define the elastic energy term
\be
f(P):=L_1 P_{ijk} P_{ijk} +L_2 P_{ikj}P_{ijk}+L_3 P_{ijj}P_{ikk}.
\ee 
In order to find relationships between  $L_1, L_2$ and $L_3$ that ensure the coercivity  we have to compute
\be
\min_{P \in \mcP, \ P_{ijk}P_{ijk}=1}f(P).
\ee
Note that we do not take into account the relation between $\mcI_2$ and $\mcI_3$ and we obtain pointwise constraints that give only {\it sufficient conditions} on the elastic constants $L_i$ to make  the energy density in $\mcF_{el}$ positive definite pointwise.

In order to derive the coercivity conditions we define the Lagrangian function
\be
F(P,\lambda):=f(P)+\lambda(P_{ijk}P_{ijk}-1)+\mu_j (P_{llj})+\nu_{ijk}(P_{ijk}-P_{jik})
\ee whose critical points satisfy the equations 
\be\label{partial:fd}
\frac{\partial F}{\partial P_{\alpha\beta\gamma}}=2(L_1+\lambda) P_{\alpha\beta\gamma}+2L_2 P_{\alpha\gamma\beta}+2L_3P_{\alpha ll}\delta_{\beta\gamma}+\mu_\gamma\delta_{\alpha\beta}+(\nu_{\alpha\beta\gamma}-
\nu_{\beta\alpha\gamma})=0,
\ee
\be\label{partial:fl}
P_{\alpha\beta\gamma}P_{\alpha\beta\gamma}=1,
\ee
\be\label{trace}
P_{\alpha\alpha\gamma}=0,
\ee
\be\label{symm}
P_{\alpha\beta\gamma}=P_{\beta\alpha\gamma}.
\ee
The last three equations are due to the constraints on $P$.

Multiplying \eqref{partial:fd} by $\delta_{\alpha\beta}$ and using \eqref{partial:fl}, \eqref{trace} and \eqref{symm} we obtain
\be
\mu_\gamma=-\frac{2}{3}(L_2+L_3)P_{\gamma ll} .
\ee
Subtracting out of \eqref{partial:fd} the same equation but with $\alpha$ and $\beta$ interchanged, while using \eqref{symm} we obtain
\be
\sigma_{\alpha\beta\gamma}:=\nu_{\alpha\beta\gamma}-\nu_{\beta\alpha\gamma}=
L_2(P_{\beta\gamma\alpha}-P_{\alpha\gamma\beta})+L_3(P_{\beta ll}\delta_{\alpha\gamma}-P_{\alpha ll}\delta_{\beta\gamma}).
\ee
Multiplying \eqref{partial:fd} respectively by $P_{\alpha\beta\gamma}$, $P_{\alpha\gamma\beta}$   and using \eqref{partial:fl} we obtain the relations
\be \label{rel1}
2(L_1+\lambda)+2L_2 P_{\alpha\gamma\beta}P_{\alpha\beta\gamma}+2L_3(P_{\alpha ll}P_{\alpha mm})=0,
\ee
\be \label{rel2}
(L_2+2(L_1+\lambda)) P_{\alpha\gamma\beta}P_{\alpha\beta\gamma}+L_2+\left(\frac{L_3-2L_2}{3}\right)(P_{\alpha ll}P_{\alpha mm})=0.
\ee
Finally, taking $\beta=\gamma$ and multiplying \eqref{partial:fd} by $P_{\alpha ll}$, while taking into account that $\delta_{\gamma\gamma}=3$, we obtain
\be \label{rel3}
\left(2(L_1+\lambda)+\frac{L_2}{3}+\frac{10}{3}L_3\right)(P_{\alpha\beta\beta}P_{\alpha ll})=0.
\ee
In order to obtain the coercivity conditions on $L_i$  it suffices that the energy computed at any critical point is positive.  We denote the elastic energy terms as
 $$
 H_1:=P_{\alpha\beta\gamma}P_{\alpha\gamma\beta}, \quad H_2:=P_{\alpha ll}P_{\alpha mm}
 $$
 and using equations \eqref{rel1}-\eqref{rel3} arrive at the following conditions:
 \begin{itemize}
\item {\it Case 1: $H_2=0$}. In this case out of \eqref{rel1} respectively \eqref{rel2} we get:
 $$
 H_1=-\frac{L_1+\lambda}{L_2}=-\frac{L_2}{L_2+2(L_1+\lambda)}
 $$ 
and  hence
 $$
 \frac{1}{H_1}=-\frac{L_2+2(L_1+\lambda)}{L_2}=-1+2H_1.
 $$ 
 Solving the quadratic equation, we have
 $$
 H_1\in\{1,-\frac{1}{2}\}.
 $$
 Using $H_1$ and $H_2$ in the energy we get the conditions:
 $$
 L_1+L_2>0, \quad 2L_1-L_2>0.
 $$
 \item {\it Case 2: $H_2\not=0$}. In this case 
$$
2(L_1+\lambda)+\frac{L_2}{3}+\frac{10}{3}L_3=0
$$
and we can find
$$
L_1+\lambda=-\frac{L_2}{6}-\frac{5}{3}L_3.
$$ 
Plugging it into \eqref{rel1}, \eqref{rel2} we obtain the following system
$$
L_2 H_1+L_3 H_2=\frac{L_2}{6}+\frac{5}{3}L_3,
$$
$$ 
\left(\frac{2L_2}{3}-\frac{10}{3}L_3\right)H_1+\left(\frac{L_3}{3}-\frac{2L_2}{3}\right)H_2=-L_2.
$$ 
Since the energy is $f(P)=L_1+L_2H_1+L_3H_2$ using the first equation we immediately obtain 
$$
L_1+\frac{L_2}{6}+\frac{5}{3}L_3>0
$$
and arrive at \eqref{sufcc}.
\end{itemize}
 It straightforward to show that the coercivity conditions in $\RR^2$ will be exactly as in \eqref{sufcc}. 

\bigskip

\noindent {\bf 2.} Now we would like to see how coercivity conditions will be changed if we use Dirichlet boundary conditions and the equivalence of elastic constants $L_2$ and $L_3$, modulo boundary terms. We define
$$
M= \frac{L_2+L_3}{2}, \ L_2 = \alpha M, \ L_3 = ({2} -\alpha ) M.
$$

Using the previous part we have:

\bea
\int_\Omega L_1Q_{ij,k}Q_{ik,j}+L_2 Q_{ik,j}Q_{ij,k}+L_3Q_{ij,j}Q_{ik,k}\,dx&=\int_\Omega L_1 Q_{ij,k}Q_{ij,k}+\underbrace{(L_2+L_3)}_{=2M} Q_{ij,j}Q_{ik,k}\,dx\non\\
&+L_2\int_\Omega Q_{ij,k}Q_{ik,j}-Q_{ij,j}Q_{ik,k}\,dx
\eea

 We note that because of Lemma~\ref{lemma:nulllagrangian} the last term is a constant, independent of $Q$, but depending just on $L_2$ and the boundary condition $Q=Q_b$ on $\partial\Omega$. Thus from the point of view of energy minimization we can ignore this constant and  therefore the elastic constants $L_2$ and $L_3$ are interchangeable. We define the set $X_\alpha$ 
$$
X_\alpha=\{ (L_1, M)\, : \, L_1 + \alpha M>0, \ 2 L_1 - \alpha M >0, \ L_1+\frac{\alpha}{6} M+\frac{5 (2-\alpha)}{3}M>0 \}
$$
and note that for each $\alpha$ the set $X_\alpha$ defines the corresponding coercivity conditions. Since $\alpha$ can be arbitrary the new coercivity conditions set $X$ is defined as $X =\cup_\alpha X_\alpha$. 
It is straightforward to obtain that 
$$
X=\{ (L_1, M)\, : \, L_1 + \frac43 M>0, L_1 >0\}.
$$
Similar considerations hold for $2D$ domains. 
\eproof

\subsection{Relation between 2D and 3D solutions}
\label{sec:2d3d}

In this section we want to investigate the relation between the critical points of the Landau-de Gennes energy in a 2D domain $\omega$ and the corresponding critical points of LDG energy in the 3D cylinder $\Omega=\omega \times (-d, d)$. It is well known that when $L_2=L_3=0$ the critical point of 2D Landau-de Gennes energy is in fact a translation invariant critical point of  the 3D Landau-de Gennes energy. Here we show that this is not necessarily the case when $L_2+ L_3 \neq 0$. 

Critical points of the 2D Landau-de Gennes energy \eqref{LDG} with Dirichlet boundary conditions 
$$
Q (x) = Q_0(x) \hbox{ for } x \in \partial \omega
$$ 
satisfy the following weak form of the Euler-Lagrange equations:
\begin{align}
 \int_{\omega} \Big[{L_1} \partial_k {Q_{ij}}  \partial_k {R_{ij}} + L_2 \partial_j Q_{ik} \partial_k R_{ij}  + L_3 \partial_k Q_{ik} \partial_j R_{ij}  + \partial f_B(Q_{ij}) R_{ij} \Big]\,dx =0,
\end{align} 
where $R(x) \in \mcS_0$ for any $x \in \omega$ and $R(x)=0$ for $x \in \partial \omega$. The summation is taken over $i,j,k=1,2,3$, and since we work in a two-dimensional domain, we assume $Q_{ij,3}=R_{ij,3}\equiv 0$.

Taking $z$-invariant boundary conditions in 3D, 
$$
Q(x,z) = Q_0(x) \hbox{ for } (x,z) \in \partial \omega \times (-d,d)
$$
and computing the first variation of the Landau-de Gennes energy we obtain
\begin{align} \label{3d}
& \int_{\Omega} \Big[{L_1} \partial_k {Q_{ij}}   \partial_k {P_{ij}} + L_2 \partial_j Q_{ik} \partial_k P_{ij}  + L_3 \partial_k Q_{ik} \partial_j P_{ij}  + \partial f_B(Q_{ij}) P_{ij} \Big]\,dx =0,
\end{align} 
where  $P(x,z) \in \mcS_0$ for any $(x,z) \in \Omega$ and $P(x,z) =0$ for any $(x,z) \in \partial \omega \times (-d,d)$.

We would like to show that  {\it unlike} in the case $L_2+L_3=0$ the minimizer in 2D will not necessarily satisfy the Euler-Lagrange equations for the 3D problem. It  is not difficult to see that \eqref{3d} produces ``natural boundary conditions"
at the top and bottom of the cylinder,
$$
\int_\omega L_1 \partial_3 Q_{ij} (x, \pm d) P_{ij}(x,\pm d) + L_2 \partial_j Q_{i3} (x,\pm d) P_{ij} (x,\pm d) + L_3 \partial_k Q_{ik} (x, \pm d) P_{i3} (x, \pm d) \, dx =0.
$$
Assuming that  the 3D minimizer $Q_{ij}$ is independent of $z$, we immediately obtain
\be\label{linearconstr}
L_2 ( \partial_j  Q_{i3}(x) + \partial_i Q_{j3}(x)) + L_3 ( \partial_k Q_{ik}(x) \delta_{3j} + \partial_k Q_{jk}(x) \delta_{3i}) =0
\ee
for all $x \in \omega$.

However, we also know that $Q$ has to satisfy Euler-Lagrange equations \eqref{eq:EL} with the boundary data $Q(x) = Q_0(x)$ for $x \in \partial \omega$. Therefore the nonlinear second order PDE has to be consistent with the linear first order PDE in \eqref{linearconstr}, and in general this is not true.  For example, it can be shown that the two-dimensional radially symmetric solutions we have considered here do not extend to $z$-independent solutions in  three dimensions if $L_2 + L_3 \neq 0$.

\subsection{Symmetry breaking}
\label{sec:symmbreaking}
In this section we compute the symmetry breaking term
\be \label{LY}
{\mathcal L} Q_{ij}=   \partial_j \partial_k Q_{ik} + \partial_i \partial_k Q_{jk}  - \frac{2}{3} \partial_l \partial_k Q_{lk} \delta_{ij}, 
\ee
where $Q$ is defined  as
$$
Q= v(r) E_0 + u(r) E_1 + w(r) E_2.
$$
It is straightforward to compute
$$
{\mathcal L} (v E_0)  = -\frac{2}{\sqrt{6}} \left(  \{Hess f\} - \frac{1}{3} \Delta v \, I_3 \right),
$$ where $f(x)=v(|x|)$.
Now we want to find ${\mathcal L} (u E_1)$ and ${\mathcal L} (u E_2)$ . We first compute the generic term (abusing notation we say $E^1=E_1$, $E^2=E_2$),
$$
\partial_j \partial_m \left( u(r) E^1_{im} \right) = (\partial_j \partial_m u(r)) E^1_{im} +  ( \partial_m u(r)) \partial_j E^1_{im} +
(\partial_j u(r)) \partial_m E^1_{im} + u(r) \partial_j \partial_m E^1_{im}.
$$
It is straightforward to obtain the following relations:
$$
\partial_j u(r) = u'(r) \frac{x_j}{r}, \quad \partial_j \partial_m u(r) = \left( u'' - \frac{u'}{r} \right) \frac{x_j x_m}{r^2} + \frac{u'}{r} \delta_{jm};\quad j,m=1,2
$$
$$
\partial_j E^1 = \frac{k x^\perp_j}{\sqrt{2} r^2} \left(  
\begin{array}{cc}
- \sin (k \varphi) & \cos(k \varphi) \\
\cos(k \varphi) & \sin(k \varphi)
\end{array}
\right), \quad
\partial_m E_{\, \cdot \, m} = \frac{k}{\sqrt{2} r} \left(
\begin{array}{c}
\cos((k-1) \varphi) \\
 \sin((k-1) \varphi)
\end{array}
\right), 
$$
where we use notation $x^\perp =(-x_2, x_1,0)$.

Using the above equalities we obtain
\bea
A_{11}=\partial_1 \partial_m \left( u(r) E^1_{1m} \right) &= \left(u'' - \frac{u'}{r} \right) \frac{x_1}{\sqrt{2} r} \cos((k-1) \varphi) +
\frac{u'}{\sqrt{2} r} \cos(k \varphi) \non \\ &+ \frac{k u'}{\sqrt{2} r} \cos((k-2) \varphi) + u \partial_1 \left(  \frac{k}{\sqrt{2} r} \cos((k-1) \varphi) \right),
\eea
\bea
A_{12}=\partial_2 \partial_m \left( u(r) E^1_{1m} \right) &= \left(u'' - \frac{u'}{r} \right) \frac{x_2}{\sqrt{2} r} \cos((k-1) \varphi) +
\frac{u'}{\sqrt{2} r} \sin(k \varphi) \non \\ &- \frac{k u'}{\sqrt{2} r} \sin((k-2) \varphi) + u \partial_2 \left(  \frac{k}{\sqrt{2} r} \cos((k-1) \varphi) \right),
\eea
\bea
A_{21}=\partial_1 \partial_m \left( u(r) E^1_{2m} \right) &= \left(u'' - \frac{u'}{r} \right) \frac{x_1}{\sqrt{2} r} \sin((k-1) \varphi) +
\frac{u'}{\sqrt{2} r} \sin(k \varphi) \non \\ &+ \frac{k u'}{\sqrt{2} r} \sin((k-2) \varphi) + u \partial_1 \left(  \frac{k}{\sqrt{2} r} \sin((k-1) \varphi) \right),
\eea
\bea
A_{22}=\partial_2 \partial_m \left( u(r) E^1_{2m} \right) &= \left(u'' - \frac{u'}{r} \right) \frac{x_2}{\sqrt{2} r} \sin((k-1) \varphi) -
\frac{u'}{\sqrt{2} r} \cos(k \varphi) \non \\ &+ \frac{k u'}{\sqrt{2} r} \cos((k-2) \varphi) + u \partial_2 \left(  \frac{k}{\sqrt{2} r} \sin((k-1) \varphi) \right).
\eea
Now we are ready to find 
\bea
{\mathcal L} (u E_1)= \left(
 \begin{array}{ccc}
A_{11}-A_{22} & A_{12} +A_{21} & 0\\
A_{12} +A_{21} & -A_{11}+A_{22}&0 \\
0 & 0& 0
\end{array}
\right) + \frac{A_{11}+A_{22}}{3} \left(
 \begin{array}{ccc}
1 & 0& 0\\
0 & 1&0 \\
0 & 0& -2
\end{array}
\right) .
\eea
Since 
$$
E_{1m}^2 =  -E_{2m}^1, \quad E_{2m}^2 = E_{1m}^1
$$
we can immediately obtain 
\bea
{\mathcal L} (u E_2)= \left(
 \begin{array}{ccc}
-A_{21}-A_{12} & -A_{22} +A_{11} & 0\\
-A_{22} +A_{11} & A_{21}+A_{12}&0 \\
0 & 0& 0
\end{array}
\right) + \frac{A_{12}-A_{21}}{3} \left(
 \begin{array}{ccc}
1 & 0& 0\\
0 & 1&0 \\
0 & 0& -2
\end{array}
\right) .
\eea

It is not difficult to compute
$$
A_{11} - A_{22} = \left( u'' + \frac{u'}{r} - \frac{k^2 u}{r^2} \right) \frac{ \cos (k\varphi)}{\sqrt{2}},
$$
$$
A_{12} + A_{21} = \left( u'' + \frac{u'}{r} - \frac{k^2 u}{r^2} \right) \frac{ \sin (k\varphi)}{\sqrt{2}},
$$
$$
A_{11} + A_{22} = \left( u'' + \frac{(2k-1) u'}{r} + \frac{k (k-2) u}{r^2} \right) \frac{ \cos ((k-2)\varphi)}{\sqrt{2}},
$$
$$
A_{12} - A_{21} = -\left( u'' + \frac{(2k-1) u'}{r} + \frac{k (k-2) u}{r^2} \right) \frac{ \sin ((k-2)\varphi)}{\sqrt{2}}.
$$
Combining the results above, we obtain
\bea
 {\mathcal L} (v(r) E_{0} ) &= \frac{1}{3} \left( v'' + \frac{v'}{r} \right) E_0 - \frac{1}{\sqrt{3}} \left( v'' - \frac{v'}{r} \right) \left( E_1 \cos((k-2) \varphi) - E_2 \sin((k-2) \varphi) \right),
\eea
\bea
{\mathcal L} (u(r) E_{1} ) &= \left( u'' + \frac{u'}{r} - \frac{k^2 u}{r^2} \right) E_1 - \frac{1}{\sqrt{3}} \left( u'' + \frac{(2k-1) u'}{r} +\frac{k (k-2) u}{r^2} \right) \cos ((k-2)\varphi)  E_0 ,
\eea
\bea
{\mathcal L} (w(r) E_{2} ) &= \left( w'' + \frac{w'}{r} - \frac{k^2 w}{r^2} \right) E_2 + \frac{1}{\sqrt{3}} \left( w'' + \frac{(2k-1) w'}{r} +\frac{k (k-2) w}{r^2} \right) \sin ((k-2)\varphi)  E_0 .
\eea

\bibliographystyle{acm}
\bibliography{LCbibliography.bib}

\end{document}